\documentclass[twocolumn,5p,sort&compress]{elsarticle}
\usepackage{amsfonts}
\usepackage{amssymb}
\usepackage{amsmath}
\usepackage{times}
\usepackage[T1]{fontenc}
\usepackage{graphicx}
\newcommand{\mi}{\mathrm{i}}
\newcommand*\diff{\mathop{}\!\mathrm{d}}
\usepackage[percent]{overpic}
\usepackage[]{matlab-prettifier}
\usepackage{filecontents}
\usepackage{dcolumn}
\newcolumntype{d}[1]{D{.}{\cdot}{#1} }
\usepackage{bm}
\usepackage{multicol}
\usepackage{relsize}
\usepackage{lipsum}

\newcommand\blfootnote[1]{%
  \begingroup
  \renewcommand\thefootnote{}\footnote{#1}%
  \addtocounter{footnote}{-1}%
  \endgroup
}
\allowdisplaybreaks
\lstMakeShortInline"

\begin{document}
\begin{frontmatter}
\title{A linear system for pipe flow stability analysis allowing for boundary condition modifications}
\author[malik]{M. Malik} \ead{dr.malik.barak@gmail.com}
\author[cf]{Martin Skote\corref{cor1}} \ead{m.skote@cranfield.ac.uk}
\cortext[cor1]{Corresponding author} 
\address[malik]{Block 32 Ghim Moh Link, \#25-298 Singapore 271032} 
\address[cf]{Cranfield University, United Kingdom} 

\begin{abstract}
An accurate system to study the stability of pipe flow that ensures regularity is presented. 
The system produces a spectrum that is as accurate as Meseguer \& Trefethen (2000), 
while providing flexibility to amend the boundary conditions without a need to modify
the formulation. The accuracy is achieved by formulating the state variables to behave as analytic
functions. We show that the resulting system retains the regular singularity at the pipe centre with 
a multiplicity of poles such that the wall boundary conditions are complemented with precisely 
the needed number of regularity conditions for obtaining unique solutions. 
In the case of axisymmetric and axially constant perturbations the computed eigenvalues match, 
to double precision accuracy, the values predicted by the analytical characteristic relations. 
The derived system is used to obtain the optimal inviscid disturbance pattern, which is found to hold 
similar structure as in plane shear flows. 
\end{abstract}

\end{frontmatter}

\section{Introduction\label{sec:intro}}
Flow through pipes is of a great importance due to vast applications ranging from cooling systems 
to fluid transportations\blfootnote{\copyright  2019. This manuscript version is made available under the CC-BY-NC-ND 4.0 license http://creativecommons.org/licenses/by-nc-nd/4.0/}. As drag reduction in such systems will enhance human life 
by lowering energy consumption, understanding 
laminar to turbulent transition has been the  subject of research since Reynolds' experiments. Despite that the flow exhibits 
instability in reality, it is linearly stable mathematically. This phenomenon has induced research to tackle the 
understanding in the perspectives of nonlinear dynamics, chaos and intermediate \textit{edge 
states}~\cite{itano2001dynamics, schneider2007turbulence, budanur2018complexity}, relatively recently.

However, for the nonlinearity to set in, the disturbances need to be brought to a 
certain finite-scale size. Algebraic transient 
growth~\cite{boberg1988onset, butler1992three, trefethen1993hydrodynamic, schmid1994optimal, 
schmid2007nonmodal} caused by 
non-normality of the linearized Navier-Stokes operator can mathematically explain this initial growth. In 
physical terms, the transient growth is described by the inviscid interaction of the mean flow and perturbation through 
the lift-up effect for the modes with azimuthal modulations and through the Orr 
mechanism of enhancing the perturbations that oppose the mean shear in the case of modes with streamwise 
modulations. Inherently, since the laminar mean flow exhibits a mean shear, the flow has a tendency to 
transfer energy into infinitesimal disturbances, which can be inflow borne or originate from 
wall roughness and vibrations. 

One of the other several routes to turbulence is distortions 
of the mean flow. 
Spatial analysis of this flow configuration predicted instability under suitable mean flow 
distortion~\cite{gavarini2004initial}. Furthermore, irregularities in the pipe geometry can serve as a cause 
for instability. For example, the flow has been found unstable in slowly diverging 
pipes~\cite{sahu2005stability}.  In addition, flow through sinusoidally corrugated pipe has exhibited linear 
instability~\cite{loh2011stability}. Wall roughness in the molecular scales has been found to be a 
destabilizing factor compared to smooth pipe, though instability has not been 
observed~\cite{pruuvsa2009influence}.

Among the early investigations that predicted the linear stability of flow through straight pipes in 
temporal setting~\cite{lessen1968stability, burridge1969comments, salwen1972stability, 
salwen1980linear}, Burridge \& Drazin~\cite{burridge1969comments} proposed an ansatz for their working 
variables for an asymptotic analysis, which was later adopted for numerical 
calculations~\cite{schmid1994optimal,o1994transient,thomas2012linear}. 

Priymak \& Miyazaki~\cite{priymak1998accurate} predicted the form of the $r-$, the radial 
coordinate, dependence of the velocity and pressure variables through an ansatz, which is suitable 
if solutions of the linear system for the perturbations that are analytic at the centreline are sought. 
It should be noted that this ansatz already incorporates the necessary behaviour for velocity fields 
deduced by Khorrami \textit{et al.}~\cite{khorrami1989application}, thus leaving only regularity as
means for further conditioning of the system.  
We would like to note that in the light of this ansatz, the 
working variables of Burridge \& Drazin~\cite{burridge1969comments}, namely, $\overline{\phi}$ 
and $\overline{\Omega}$ (defined in section~\ref{subsec:validations}), which are related to the 
radial velocity and radial vorticity have the property, $(\overline{\phi}, \overline{\Omega})\sim (r^2, 
\textrm{const})$ for $n = 0$ to the leading order close to the centreline, and  
$(\overline{\phi}, \overline{\Omega})\sim (r^{|n|}, r^{|n|})$ for $n\neq 0$, where $n$ is the 
azimuthal wavenumber. This suggests that for higher values of $|n|$ (e.g., $|n| \geq 5$), all 
the eigenfunctions will be vanishingly small in a finite neighbourhood of the centreline. 
Consequently, there will be precision loss in determining the eigenfunctions. This complication 
arises because the eigenvalues are multiplied by exceedingly low values in the 
discretized eigensystem equations close to the the centreline.

Furthermore, an increase of collocation points would only worsen the situation by 
enhancing the round-off errors since the eigenfunctions $(\overline{\phi}, \overline{\Omega})^{\rm T}$ do 
not undergo steep variations close to the centreline. These round-off errors at the end of the iterative 
procedures for eigenvalues demonstrate themselves as pseudospectra. Since pseudospectra generally affects modes 
decaying at higher rate than the least decaying ones, the optimal patterns can get distorted at such 
high values of $n$.

Priymak \& Miyazaki~\cite{priymak1998accurate} also developed a method for solving the linearized 
Navier-Stokes equations cast as a system of two unknowns, the radial and azimuthal perturbation 
velocities, that utilizes the above findings, thus producing very accurate spectra. However, the 
linear operator itself needs to be determined in a complex algorithmic way, unlike the case of Burridge 
\& Drazin~\cite{burridge1969comments}. In~\cite{priymak1998accurate}, the linear system is determined 
through a sequence of steps numerically, where the pressure needs to be found by inverting a Poisson 
operator, and the continuity is being imposed by a correcting pressure which leads to correcting 
velocities, similar to the steps of the SIMPLE algorithm.

Meseguer \& Trefethen~\cite{meseguer2000spectral,meseguer2003linearized} deployed the ansatz of 
Priymak and Miyazaki, 
and have produced very accurate results with a number of collocation points as low as $55$. 
Their work is in Petrov-Galerkin formalism which uses solenoidal trial basis and test functions 
that satisfies the particular boundary condition of no-slip. However, this method comes with a 
task of changing the trial and test bases if the boundary conditions are other than no-slip. 
In some occasions, the boundary conditions are even time dependent (see for 
example,~\cite{malik2018growth}). Finding appropriate solenoidal trial and test bases that 
satisfy arbitrary boundary conditions and the properties prescribed by the ansatz of Priymak \& 
Miyazaki is not a simple task, in general. 

The above facts motivate us to formulate a method that is as accurate and efficient as Priymak \& 
Miyazaki~\cite{priymak1998accurate} and Meseguer \& Trefethen~\cite{meseguer2000spectral}, but as 
flexible and simple as a usual spectral collocation method with 2-tuple state variable similar to that 
of Ref.~\cite{schmid1994optimal} or~\cite{burridge1969comments}. To this end, we use the ansatz of 
Priymak \& Miyazaki for velocity fields to determine a 2-tuple working variables that do not go like a 
power law close to the centreline. As the boundary conditions will be imposed on the unknowns unlike the 
Petrov-Galerkin method, 
the formulated 
equations are applicable for a range of boundary conditions.

The theoretical derivations are presented in section~\ref{sec:linsys}, while the numerical results
are discussed in section~\ref{sec:numerical}.
Finally, in section~\ref{sec:ellingsen}, we derive the Ellingsen and Palm 
solutions~\cite{ellingsen1975stability} that are valid for inviscid axially constant 
modes that captures nonmodal algebraic growth. Derivation of this is remarkably simpler 
in the working variables in this paper. Such solutions have been found useful in the 
case of plane-shear flows to compute the optimal patterns that demonstrate the lift-up effect 
leading to streaks. We find that the pipe flow configuration retains these features. Although 
this is  known through earlier viscous computations~\cite{schmid1994optimal} and 
DNS~\cite{zikanov1996instability}, we arrive at these results in the inviscid limit.

\section{Theoretical derivations\label{sec:linsys}}
Let us consider the cylindrical coordinates, $\boldsymbol{r} = (x,r,\theta)$, 
which are the axial, radial and azimuthal coordinates, respectively. We emphasize the unusual order of the coordinates used here. Let $\boldsymbol{u}(\boldsymbol{r}) = (u, v, w)^{\rm T}$ be the velocities normalized with respect to the centreline velocity. Let us use the decomposition, 
$\boldsymbol{u} = \boldsymbol{U} + \boldsymbol{\tilde{u}}$ where the mean flow, 
$\boldsymbol{U} = (U, 0, 0)^{\rm T}$ and $U = 1-r^2$, and the considered perturbation, 
$\boldsymbol{\tilde{u}}$ is such that $\Vert\boldsymbol{\tilde{u}}\Vert \ll \Vert\boldsymbol{U}\Vert$, which are governed by,
\begin{equation}
\boldsymbol{\tilde{u}}_t +  (\boldsymbol{U}\cdot\boldsymbol{\nabla}) \boldsymbol{\tilde{u}} +  (\boldsymbol{\tilde{u}}\cdot\boldsymbol{\nabla}) \boldsymbol{U} = -\boldsymbol{\nabla}\tilde{p} + Re^{-1}\nabla^2\boldsymbol{\tilde{u}}, \label{eq:pert}
\end{equation}
and by the continuity condition, $\boldsymbol{\nabla}\cdot\boldsymbol{\tilde{u}} = 0$. In 
Eq.~(\ref{eq:pert}), $Re = U_cR/\nu$, where $U_c$ is the dimensional centreline velocity, $R$ is 
the pipe radius, and $\nu$ is the kinematic viscosity. The $\tilde{p}$ in Eq.~(\ref{eq:pert}) is the 
perturbation pressure. Substituting $(\boldsymbol{\tilde{u}}, \tilde{p}) = (\boldsymbol{u'}, p') \exp[\mi(\alpha x + n\theta - \omega t)]$ into Eq.~(\ref{eq:pert}) to consider a single Fourier mode, 
and upon taking the curl to remove the appearance of pressure, we obtain the following three equations:
\begin{multline}
\!\mi (\alpha U-\omega)(Du'-\mi\alpha v') = -\mi \alpha U_ru' - (U_{rr}+U_rD)v' \\+ Re^{-1}\left ( \Delta (Du'-\mi\alpha v') -2\mi n r^{-2}\eta' \right ) \label{eq:curl:1}
\end{multline}
\begin{multline}
\! \mi(\alpha U-\omega)[D(rw')-\mi nv'] = -\mi\alpha rU_rw' \\+  Re^{-1}\left ( \overline{\Delta} [D(rw') - \mi nv'] \right ) \label{eq:curl:2}
\end{multline}
\begin{multline}
\!\mi(\alpha U-\omega)\eta' = -\mi nU_r r^{-1}v' \\+ Re^{-1} \left ( \Delta\eta' + 2 n r^{-2} (\alpha v' + \mi Du') \right ) \label{eq:curl:3}
\end{multline}
where $\Delta = \left[D^2 + r^{-1}D -r^{-2}(d+1)\right]$, $\overline{\Delta} = \left[D^2 - r^{-1}D -r^{-2}(d-1)\right]$ , $\eta' = \mi(nr^{-1}u'-\alpha w')$ and $d = n^2 + \alpha^2r^2$. As a matter of 
convention, we use suffix $\scriptsize{r}$ to imply differentiation of mean flow variables with respect 
to $r$, and the operator $D$ to represent the same for the fluctuation variables. Using the above 
definition of $\eta'$ and the continuity equation, the variables $u'$ and $w'$ can be written as 
\begin{align}
u' =& \mi r d^{-1}(\alpha v' + \alpha r Dv' -n\eta') \quad \text{and} \label{eq:uprime}\\
w' =& \mi d^{-1}(n v' + n r Dv' +\alpha r^2\eta'). \label{eq:wprime}
\end{align}
Priymak \& Miyazaki~\cite{priymak1998accurate} observed the following behaviour close to the centreline.
\begin{equation}
(u', v', w') = \left \{
\begin{array}{lll}
(\psi_1, r\phi, r\psi_2) &\text{for}  &n = 0,\\
(r^{\ell+1}\psi_1,r^\ell\phi,r^\ell\psi_2) &\text{for}  &n \neq 0
\end{array}
\right . \label{eq:uvwprimes_r_variations}
\end{equation}
where $\ell = |n|-1$ and, $\psi_1$, $\psi_2$ and $\phi$ are analytic functions having Taylor 
expansion around the centreline with vanishing coefficients of odd powers of $r$. 
Substituting Eq.~(\ref{eq:uvwprimes_r_variations}) into the definition of $\eta'$, we get its 
behaviour. In summary, $v'$ and $\eta'$ have the following forms:
\begin{equation}
(v', \eta') = \left \{
\begin{array}{lll}
(r\phi, r\Omega) &\text{for}  &n = 0 \quad \text{and}\\
(r^\ell\phi,r^\ell\Omega) &\text{for}  &n \neq 0
\end{array}
\right . \label{eq:vetaprimes_r_variations}
\end{equation}
with $\Omega$ an analytical function having similar Taylor expansion as $\phi$. Substituting 
Eq.~(\ref{eq:vetaprimes_r_variations}) into Eq.~(\ref{eq:uprime}) and Eq.~(\ref{eq:wprime}), we get 
\begin{align}
u' =&  \left \{
\begin{array}{ll}
\frac{\mi r^{\ell+1}}{d}[\alpha (\ell+1)\phi + \alpha r D\phi -n\Omega] &n \neq 0 \\
\frac{\mi}{\alpha}(2\phi + r D\phi) &n = 0
\end{array}  \label{eq:uprime_new}
\right .	\\
w' =& \left \{
\begin{array}{lll}
\mi r^\ell d^{-1}[n (\ell + 1)\phi + n r D\phi +\alpha r^2\Omega] &n \neq 0,\\
\mi  \alpha^{-1}r\Omega &n = 0.
\end{array}
\right .	 \label{eq:wprime_new}
\end{align}
First, let us consider the case of $n \neq 0$. Upon substituting Eqs.~(\ref{eq:vetaprimes_r_variations})-(\ref{eq:wprime_new}) into Eqs.~(\ref{eq:curl:1})-(\ref{eq:curl:3}) we get
\begin{multline}
(\omega - \alpha U) \left ( \alpha L_1\phi -nL_2\Omega \right ) = \alpha L_3\phi  -\alpha nL_4\Omega \\ + \mi Re^{-1} \left (\alpha L_5\phi  -nL_6\Omega\right ), \label{eq:phi_from_eq_uprime}
\end{multline}
\begin{multline}
(\omega - \alpha U) \left ( n L_7\phi + \alpha L_8\Omega \right ) = n L_9\phi  -\alpha^2 L_{10}\Omega \\ + \mi Re^{-1} \left (n L_{11}\phi  +\alpha L_{12}\Omega\right ), \quad \text{and} \label{eq:phi_from_eq_wprime}
\end{multline}
\begin{equation}
(\omega - \alpha U)\Omega = nr^{-1}U_r\phi + \mi Re^{-1} \left ( 2\alpha n L_{13} \phi + L_{14}\Omega \right ) \label{eq:omega}
\end{equation}
where the operators, $L_{1 - 14}$ are given in~\ref{app:oper}. Eq.~(\ref{eq:phi_from_eq_uprime}) and 
Eq.~(\ref{eq:phi_from_eq_wprime}) are fourth order in $\phi$ and third order in $\Omega$. As we have 
three equations, Eqs.~(\ref{eq:phi_from_eq_uprime})-(\ref{eq:omega}) for two unknowns, $\phi$ and 
$\Omega$, we use Eq.~(\ref{eq:phi_from_eq_uprime}) and Eq.~(\ref{eq:phi_from_eq_wprime}) to reduce the 
order of $\Omega$. The resulting equation that is fourth order in $\phi$ and second order in $\Omega$ is
\begin{multline}
(\omega - \alpha U) \left ( rdL_1\phi + 2\alpha n r^2d^{-1}\Omega \right ) = \alpha r (U_r-rU_{rr})\phi \\+ \mi Re^{-1} \left (L_{15}\phi  +4\alpha nL_{16}\Omega\right ), \label{eq:phi} 
\end{multline}
where $L_{15}$ and $L_{16}$ are operators given in~\ref{app:oper}. 

Four more tasks are carried out before arriving at the required final system to work with: 
Firstly, the appearance of second derivative of $\Omega$ in Eq.~(\ref{eq:phi}) is removed with 
the help of Eq.~(\ref{eq:omega}). Though this will not reduce the order of the system, it helps 
in deriving simpler conditions of regularity which are shown later. The outcome is that the 
order of the derivatives of $\Omega$ in the yet-to-be-derived regularity condition is less than 
the order of the governing equation. 

Secondly, the even parity nature of $\phi$ and $\Omega$, which have Taylor series expansions 
$\phi(r) = \sum_n a_n r^{2n}$ and $\Omega(r) = \sum_n b_n r^{2n}$ is implemented by the 
transformation, $y = r^2$. Since this will transform $\phi$ and $\Omega$ to general analytic 
functions, Chebyshev polynomials of both odd and even orders can be used for the spectral expansion 
without any loss of efficiency. (For another method, where only even Chebyshev polynomials 
are used, see Priymak \& Miyazaki~\cite{priymak1998accurate} and Meseguer \& 
Trefethen~\cite{meseguer2000spectral}). 

Thirdly, similar operations performed for the case of $\{n=0,\alpha \neq 0\}$ by substituting 
$\phi$ and $\Omega$ from Eq.~(\ref{eq:vetaprimes_r_variations}), give the required system for 
axisymmetric disturbances. However, it should be noted that this system is derived without using 
Eq.~(\ref{eq:curl:2}) as there will not be a need to reduce the order of $\Omega$ in the system. 
The obtained system from Eq.~(\ref{eq:curl:1}) and Eq.~(\ref{eq:curl:3}) will already be fully 
decoupled in this case with six as the order of the system.

Finally, a system suitable for the case of $\{n=0, \alpha = 0\}$, the axisymmetric and axially constant modes, is derived with $y$ $(\equiv r^2)$ as the 
independent variable. This special case requires a different set of working variables, as 
$v'(y) = 0$ (hence, $\phi = 0$) for all $y$ due to continuity, and $\eta'(y) = 0$ (hence, 
$\Omega = 0$) by definition. We choose the variables, $\psi_1$ and $\psi_2$ defined in 
Eq.~(\ref{eq:uvwprimes_r_variations}) for this purpose. The required system for this case can be 
derived from the axial and azimuthal components of  Eq.~(\ref{eq:pert}).

\subsection{Linear System \label{linsys_summary}}
In summary, the final system of equations are
\begin{multline}
(\omega -\alpha U) \left ( [-\alpha^2 g_1 + g_2 \mathcal{D} + 4yd^3\mathcal{D}^2]\phi -2\alpha n d^2\Omega \right ) \\ =-8\alpha d^2n^2 U_y\phi + \mi Re^{-1} \left ([\alpha^4 g_3 -8\alpha^2 g_4 \mathcal{D}  \right . \\ \left . + g_5 \mathcal{D}^2+8y(g_2+4d^3)\mathcal{D}^3 + 16y^2d^3\mathcal{D}^4]\phi  \right . \\ \left .  - 8\alpha^3 n [d_2 + 2dy\mathcal{D}]\Omega \right ), \label{eq:phi_in_y}
\end{multline}
\begin{multline}
(\omega - \alpha U)d^2\Omega = 2nd^2U_y\phi + \mi Re^{-1}\left ( 2\alpha n\{\alpha^2g_6 \right . \\ \left . -2[d_2 + (\ell+3)d]\mathcal{D} -4yd\mathcal{D}^2\}\phi  +\{-\alpha^2g_7 \right . \\  \left . + 4d[(\ell+1)d+n^2]\mathcal{D} + 4yd^2\mathcal{D}^2 \} \Omega  \right ), \label{eq:omega_in_y}
\end{multline}
for the case of $n\neq0$, and 
\begin{multline}
(\omega-\alpha U)(8\mathcal{D} +4y\mathcal{D}^2-\alpha^2)\phi = \mi Re^{-1} \left ( 16y^2\mathcal{D}^4 \right . \\ \left . +96y\mathcal{D}^3+(96-8\alpha^2y)\mathcal{D}^2-16\alpha^2\mathcal{D} +\alpha^4 \right )\phi \label{eq:phi_in_y_n_0}
\end{multline}
\begin{equation}
(\omega-\alpha U)\Omega = \mi Re^{-1}\left ( 4y\mathcal{D}^2 + 8\mathcal{D} -\alpha^2 \right )\Omega \label{eq:omega_in_y_n_0}
\end{equation}
for the case of $\{n = 0,\alpha \neq 0\}$, and
\begin{align}
-\mi\omega\psi_1 =& 4 Re^{-1} \left ( y\mathcal{D}^2+\mathcal{D}\right )\psi_1 \label{eq:psi1_in_y}\\
-\mi\omega\psi_2 =& 4 Re^{-1} \left ( y\mathcal{D}^2+2\mathcal{D}\right )\psi_2 \label{eq:psi2_in_y}
\end{align}
for the case of $\{n = 0,\alpha = 0\}$, where $y \equiv r^2$, 
$\mathcal{D}=\frac{\mathrm{d}}{\mathrm{d}y}$ and the $d_2$ and $g_i(y)$ ($i = 1\cdots7$) are 
functions as given in~\ref{app:oper}. The 
order of the system is six for each of the cases of $n\neq 0$ and $\{n=0,\alpha\neq0\}$, whereas 
it is four in the case of $\{n=0,\alpha=0\}$. In the latter case, the order is only four owing to 
the fact that we did not have to take the curl of Eq.~(\ref{eq:pert}) to remove the appearance of 
$\boldsymbol{\nabla} p'$, as pressure is a constant in the axial and azimuthal 
directions similar to the other perturbation quantities. 
All modes of this flow configuration, for 
all cases of $\alpha$ and $n$, are observed to be decaying. For the special case of $\{n=0,\alpha=0\}$, a physical reason can be deduced since the time evolution is dictated only by the viscous 
dissipation because both the mean shear and the pressure gradient terms are zero. 

\subsection{Boundary and Regularity Conditions}
The systems given by Eqs.~(\ref{eq:phi_in_y})-(\ref{eq:omega_in_y}) and 
Eqs.~(\ref{eq:phi_in_y_n_0})-(\ref{eq:omega_in_y_n_0}) need to be solved with six conditions for 
a unique solution, and four conditions for the system of 
Eqs.~(\ref{eq:psi1_in_y})-(\ref{eq:psi2_in_y}). 
In the case of $n \neq 0$ and $\{n = 0,\alpha \neq 0\}$, three of the needed six conditions are 
straight forward: they are the no-slip and no-penetration conditions at the pipe boundary. These 
conditions in terms of the variables $\phi$ and $\Omega$ read as 
\begin{equation}
\phi = \mathcal{D}\phi = \Omega = 0 \quad \text{at} \quad y = 1. \label{eq:boundaryconditions}
\end{equation}
The other three conditions should come from the regularity of solutions at $y = 0$. Note that, 
although we are seeking analytic solutions by making the required transformations, 
$\{v', \eta'\} \rightarrow \{\phi, \Omega\} $ with relations given in 
Eq.~(\ref{eq:vetaprimes_r_variations}), the appearance of regular singularity in 
Eq.~(\ref{eq:phi_in_y})-(\ref{eq:omega_in_y_n_0}) implies that they still can allow 
non-analytic solutions such as, for analogy, the Bessel functions of second kind $Y_\nu$ for 
positive integer $\nu$ in the case of the Bessel equation. In a numerical procedure, the 
transformation in Eq.~(\ref{eq:vetaprimes_r_variations}) and the transformation of the coordinate, 
$y=r^2$ would only guarantee the accuracy by shifting the focus to the factors of $v'$ and $\eta'$, 
namely, $\phi$ and $\Omega$, whose $r-$variations are unknown, and allow us to work with  a reduced number 
of Chebyshev polynomials in the spectral expansion, while still allowing non-analytic solutions. 
Therefore, we need regularity conditions to explicitly rule out non-analytic solutions in our solution 
procedure. 

Since we seek $\phi$ and $\Omega$ as analytic functions, they have Taylor series expansion around $y=0$. 
This warrants all orders of the derivatives to be of $O(1)$, which implies that all derivative terms in 
Eq.~(\ref{eq:phi_in_y})-(\ref{eq:omega_in_y_n_0}) that are multiplied by $y$ should vanish in the limit 
$y \rightarrow 0$. For analogy, such requirement is satisfied by the the analytic $J_\nu$, the Bessel 
functions of first kind for positive integer $\nu$, in the case of the Bessel equation. This instantly 
gives us two conditions of regularity for each of the cases, which are
\begin{multline}
(\omega -\alpha U) \left ( [-\alpha^2 g_1 + g_2 \mathcal{D}]\phi -2\alpha n d^2\Omega \right ) = \\ -8\alpha n^2d^2 U_y \phi + \mi Re^{-1} \left ([\alpha^4 g_3 -8\alpha^2 g_4 \mathcal{D} \right . \\ \left . + g_5 \mathcal{D}^2]\phi - 8\alpha^3nd_2\Omega \right ), \label{eq:phi_in_y_regularity}
\end{multline}
\begin{multline}
(\omega - \alpha U)d^2\Omega = 2nd^2U_y\phi + \mi Re^{-1}\left ( 2\alpha n\{\alpha^2g_6 \right . \\ \left . -2[d_2 + (\ell+3)d]\mathcal{D}\}\phi +\{-\alpha^2g_7 \right . \\ \left . + 4d[(\ell+1)d+n^2]\mathcal{D} \} \Omega  \right ), \label{eq:omega_in_y_regularity}
\end{multline}
for the case of $n\neq0$, and 
\begin{multline}
(\omega-\alpha U)(8\mathcal{D} -\alpha^2)\phi = \mi Re^{-1} \times \\ \left ( 96\mathcal{D}^2-16\alpha^2\mathcal{D} +\alpha^4 \right )\phi \label{eq:phi_in_y_n_0_regularity}
\end{multline}
\begin{equation}
(\omega-\alpha U)\Omega = \mi Re^{-1}\left (8\mathcal{D} -\alpha^2 \right )\Omega \label{eq:omega_in_y_n_0_regularity}
\end{equation}
for the case of $\{n=0, \alpha\neq0\}$. We need one more condition for each of these cases of 
$n\neq 0$ and $\{n=0, \alpha\neq0\}$. Note that in Eq.~(\ref{eq:phi_in_y}) 
and~(\ref{eq:phi_in_y_n_0}) the term $\mathcal{D}^4\phi$ is being multiplied by $y^2$. Requiring 
$\lim_{y\rightarrow 0}y^2\mathcal{D}^4\phi = 0$ would only yield the constraint that 
$y\mathcal{D}^4\phi \sim O(1)$. Therefore, the other condition is the one that would enforce 
$\lim_{y\rightarrow 0}y\mathcal{D}^4\phi = 0$, since $\phi$ is analytic. The following two facts 
should be noted in order to derive the required condition: (1) As the Eq.~(\ref{eq:phi_in_y}) 
and~(\ref{eq:phi_in_y_n_0}) are both valid for the whole domain, their derivatives with respect 
to $y$ are also valid through out the flow field; (2) Due to regular singularity, 
Eq.~(\ref{eq:phi_in_y_regularity}) and~(\ref{eq:phi_in_y_n_0_regularity}) are at most second 
order in $\phi$ and zero-th order in $\Omega$, allowing them to be treated as boundary conditions. 
Furthermore, the derivatives of Eq.~(\ref{eq:phi_in_y}) and~(\ref{eq:phi_in_y_n_0}) with 
respect to $y$ at the centreline can be used as boundary conditions since they are of at most first 
order in $\Omega$ and third order in $\phi$, which are less by one order comparing to the orders of these variables in the governing equations. 
Hence, the additional regularity condition that guarantees the differentiability of 
$\phi$ until fourth order for each cases of $n\neq0$ and $\{n=0,\alpha \neq 0\}$ are the 
following:
\begin{multline}
(\omega-\alpha)\left ( [-\alpha^4g_8 + \alpha^2g_9\mathcal{D} + g_{10}\mathcal{D}^2]\phi -2\alpha n^3[2\alpha^2 \right . \\ \left . + n^2\mathcal{D}]\Omega \right ) = [\alpha^3g_{11} -4\alpha n^6\ell\mathcal{D}]\phi  +2\alpha^2n^5(1-D)\Omega \\ + \mi Re^{-1} \left ( [\alpha^6g_{12} +\alpha^4 g_{13}\mathcal{D} -\alpha^2g_{14}\mathcal{D}^2 + g_{15}\mathcal{D}^3]\phi \right . \\ \left . -8\alpha^3n[\alpha^2(\ell-1)+n^2(\ell+3)\mathcal{D}]\Omega\right ) \label{eq:phi_in_y_derivative_regularity}
\end{multline} 
for $n\neq0$, and
\begin{multline}
(\omega-\alpha)(12\mathcal{D}^2 - \alpha^2\mathcal{D})\phi = \alpha(\alpha^2-8\mathcal{D})\phi \\+ \mi Re^{-1} (192\mathcal{D}^3 -24\alpha^2\mathcal{D}^2 +\alpha^4\mathcal{D})\phi \label{eq:omega_in_y_derivative_regularity}
\end{multline}
for the case of $\{n=0, \alpha\neq0\}$, where $g_{8\cdots15}$  are constants defined in~\ref{app:oper}. In 
Eq.~(\ref{eq:phi_in_y_derivative_regularity}), we have substituted the values of mean flow variables for 
brevity.

Note that the multiplicity of the pole in the governing equation is such 
that it yields exactly three conditions to complement the three boundary conditions at the wall.
Had we instead derived another regularity condition by differentiating the governing equations once more 
and evaluated at the pipe centre, the condition would not have served as a boundary condition 
since it would have had the
same orders of derivatives of the unknowns as they appear in the governing equations. 
Similarly, had the 
multiplicity been lower than in the present case, we would not have been able to obtain a total of six 
conditions to solve the sixth-order system.

In the case of $\{n=0, \alpha=0\}$, the boundary and regularity conditions are given by 
\begin{equation}
\psi_1 = \psi_2 = 0 \quad \text{at} \quad y=1, \quad \text{and},  \label{eq:psi_boundary}
\end{equation}
\begin{equation}
\left ( \omega-4\mi Re^{-1} \mathcal{D} \right )\psi_1 = \left ( \omega-8\mi Re^{-1} \mathcal{D} \right )\psi_2 = 0. \label{eq:psi_regularity}
\end{equation}
at $y = 0$.The regularity condition for $\psi_2$ in Eq.~(\ref{eq:psi_regularity}) and the governing 
equation Eq.~(\ref{eq:psi2_in_y}) can also be obtained from the limit 
$\alpha\rightarrow 0$ in Eq.~(\ref{eq:omega_in_y_n_0_regularity}) and 
Eq.~(\ref{eq:omega_in_y_n_0}), respectively, since $\psi_2 = \mi \alpha^{-1}\Omega$ for $n=0$. 
Therefore, one can expect that a part of the spectrum for the case of $\{n=0, \alpha = 0\}$ can 
be obtained from the system for the case of $\{n=0, \alpha \neq 0\}$ in the limit of 
$\alpha\rightarrow 0$. However, there are a new set of modes for this case that cannot be 
obtained from the system for finite $\alpha$ in the limit $\alpha\rightarrow 0$. These modes 
originate from Eq.~(\ref{eq:psi1_in_y}). Both sets of these modes are governed by Sturm and Liouville 
theory (see for example~\cite{boyce2009elementary}) since the underlying operators are self-adjoint. In 
the next subsection, the solutions for this case $\{\alpha = n= 0\}$ are derived in $y-$ coordinate so 
that the characteristic relation are presented in a form that is less susceptible to numerical errors. 
This will facilitate to validate the eigenvalues obtained from the numerical spectral collocation solutions
of Eq.~(\ref{eq:psi1_in_y}) and~(\ref{eq:psi2_in_y}). 

It should be highlighted that the boundary conditions can be changed depending on the flow 
configurations, without affecting the formulation, i.e., Eq.~(\ref{eq:phi_in_y})-(\ref{eq:psi2_in_y}). 
This is the advantage of the present system compared with Meseguer and 
Trefethen~\cite{meseguer2000spectral}, where one would be tasked with finding appropriate basis and test 
functions that satisfies the boundary and regularity conditions, besides satisfying the other 
condition of continuity.

\subsection{Stokes Modes with no Transient Growth}
The operators on the RHS in Eq.~(\ref{eq:psi1_in_y})-~(\ref{eq:psi2_in_y}) are examples of Stokes 
operators. The general Stokes operators are defined as Laplacians acted upon by Leray projectors which 
ensures continuity~\cite{foias2001navier}. In the present case of $\{\alpha = 0, n=0\}$, the continuity 
is automatically satisfied without any constraint as $v' = 0$, and hence $\boldsymbol{\nabla} p' =0$. 
This implies that the Leray projector is an identity operator in the present case. It is well known 
that the eigenvalues of a general Stokes operators are decaying and that the operator itself is 
self-adjoint for no-slip boundary conditions. For this particular case of $\{ \alpha = 0, n= 0 \}$, 
the solutions and characteristic relation are known to be Bessel functions of 
$r\sqrt{|\Im\{\omega\}|Re}$ and their roots at $r =1$~\cite{salwen1972stability,rummler1997eigenfunctions}. (For such modes of plane Poiseuille flow, 
see~\cite{rummler1997eigenfunctionsii}. For application of these Stokes modes, 
see~\cite{batcho1994generalized}.)

Eq.~(\ref{eq:psi1_in_y})-~(\ref{eq:psi2_in_y}) can be rewritten as
\begin{align}
\lambda\psi_1 =& \left ( y\mathcal{D}^2+\mathcal{D}\right )\psi_1, \label{eq:psi1_in_y_lambda}\\
\lambda\psi_2 =& \left ( y\mathcal{D}^2+2\mathcal{D}\right )\psi_2, \label{eq:psi2_in_y_lambda}
\end{align}
where $\lambda = -\mi \omega Re/4$. The differential operators on the RHS are self-adjoint under the inner product defined as the cross-sectional area integral of velocity fields, i.e., $\psi_1$ or $\sqrt{y}\psi_2$. The Sturm and Liouville theory suggests that the eigenvalues, $\lambda_j$ of Eq.~(\ref{eq:psi1_in_y_lambda}) and $\lambda_k$ of Eq.~(\ref{eq:psi2_in_y_lambda}) with $j,k \in 1, 2, 3, \cdots$ are real and the eigenfunctions, $\psi_{1,j}$ and $\psi_{2,j}$ are orthogonal in the sense,
\begin{equation}
\left . 
\begin{array}{l}
\int_0^1 \psi_{1,j}^* \psi_{1,k}^{} \diff y = C_1 \ \delta_{jk}\\
\int_0^1 \psi_{2,j}^* \psi_{2,k}^{} y \diff y = C_2 \ \delta_{jk}
\end{array}
\right \}
\label{eq:orth}
\end{equation} 
where $\delta_{jk}$ is the Kronecker delta and, $C_1$ and $C_2$ are constants.

As we are interested in solutions of Eq.~(\ref{eq:psi1_in_y_lambda})-(\ref{eq:psi2_in_y_lambda}) that 
are analytic at the centreline, we can bypass the method of Frobenius and resort to a simple power series 
method. Let $\psi_1 = \sum_{k\geq0}\overline{a}_ky^k$ and $\psi_2 = \sum_{k\geq0}\overline{b}_ky^k$, 
where the constants, $\overline{a}_k$ and $\overline{b}_k$, are to be determined. Substituting these 
into Eq.~(\ref{eq:psi1_in_y_lambda})-(\ref{eq:psi2_in_y_lambda}), and matching the coefficients of like 
powers of $y$, we get the eigenfunctions, $\psi_{1,l}$ and $\psi_{2,l}$ as,
\begin{align}
\psi_{1,l}^{}(y) =& \overline{a}_0\left ( 1 + \sum_{k\geq1} \left ( \lambda_{1,l}^{}\,y\right )^k (k!)^{-2} \right ) \quad \text{and} \label{eq:psi1_soln}\\
\psi_{2,l}^{}(y) =& \overline{b}_0\left (1+\sum_{k\geq1} \left ( \lambda_{2,l}^{}\,y \right )^k [k!(k+1)!]^{-1} \right ), \label{eq:psi2_soln}
\end{align}
where $\lambda_{1,l}$ with $l = 1,2,3, \cdots$ are the real solutions of the characteristic equation, 
\begin{equation}
1+\sum_{k\geq1}  \lambda^k  (k!)^{-2} = 0, \label{eq:psi1_poly}
\end{equation}
and $\lambda_{2,l}$ are the real solutions of 
\begin{equation}
1+\sum_{k\geq1}  \lambda^k  [k!(k+1)!]^{-1} = 0, \label{eq:psi2_poly}
\end{equation}
which arise due to the no-slip conditions given by Eq.~(\ref{eq:psi_boundary}). The 
Eqs.~(\ref{eq:psi1_poly}) and~(\ref{eq:psi2_poly}) can also be written as the roots of 
$J_0(\sqrt{|\Im\{\omega\}|Re})=0$ and 
$J_1(\sqrt{|\Im\{\omega\}|Re})=0$~\cite{salwen1972stability,rummler1997eigenfunctions}, however the 
above expressions are helpful to obtain the eigenvalues accurately as the roots of polynomials that 
approximate the LHS of Eqs.~(\ref{eq:psi1_poly}) and~(\ref{eq:psi2_poly}) by cutting off the summation at 
some $k=k_{\max}$. The constants, $\overline{a}_0$ and $\overline{b}_0$ in Eqs.~(\ref{eq:psi1_soln}) 
and~(\ref{eq:psi2_soln}), respectively, can be of any value as long as the linearity of the perturbations 
with respect to the mean flow is respected. Without loss of generality they can be the normalization 
constants, $\overline{a}_0 = 1/\sqrt{(\psi_1,\psi_1)}$ and $\overline{b}_0 = 1/\sqrt{\langle\psi_2,\psi_2\rangle}$, 
where the inner-products are defined as 
$(\psi_1,\ \psi_1) = \int_0^1 \psi_1^*(y)\psi_1^{}(y)\diff y$ and 
$\langle\psi_2,\ \psi_2\rangle = \int_0^1 y\psi_2^*(y)\psi_2^{}(y)\diff y$. It should be mentioned that 
$\int_0^1 \diff y$ is in fact an area element, $2\int_0^1 r \diff r$. Under such normalizations, the constants $C_1$ and $C_2$ of Eq.~(\ref{eq:orth}) become, $C_1 = C_2 = 1$. The convergence of the series
in Eq.~(\ref{eq:psi1_soln}) and (\ref{eq:psi2_soln}) is obvious by comparison tests with the series for 
$\exp[ \lambda_{1,l}^{}\,y]$ and $\exp[ \lambda_{2,l}^{}\,y]$, respectively. The implication of the 
orthogonality is that superposition of such modes cannot have transient growth. A velocity state can be 
represented using the notation introduced at the beginning of this section as $\boldsymbol{u'}_{kl}(y) = \{\psi_{1,k}(y), 0, \sqrt{y}\psi_{2,l}(y)\}^{\mathrm{T}}$. Then the superposition velocity is given by,
\begin{equation}
\boldsymbol{\tilde{u}}(t,y) = \exp(\mi[\alpha x + n\theta])\sum_{k,l} \boldsymbol{a}_{kl} \exp{(\boldsymbol{\Lambda}_{kl}t)}\boldsymbol{u'}_{kl},
\end{equation}
where $\boldsymbol{\Lambda}_{kl} = \mathrm{diag}\{\lambda_k, 0 , \lambda_l\}$ and $\boldsymbol{a}_{kl} = \mathrm{diag}\{a_k, 0 , a_l\}$ are diagonal matrices for each chosen $k$ and $l$, where in turn 
$\lambda_k$ and $\lambda_l$ are the real solutions of Eq.~(\ref{eq:psi1_poly}) 
and~(\ref{eq:psi2_poly}), respectively, and $a_k$ and $a_l$ are some constants.

Let us define the energy 
as, $E(t) =2^{-1}\int_0^1 |\boldsymbol{\tilde{u}}|^2dy$.
\begin{multline}
E(t) = 2^{-1}\sum_{k,l} \sum_{i,j} \boldsymbol{a}_{ij}^\text{H} \boldsymbol{a}_{kl}^{} \\ \times \exp{([\boldsymbol{\Lambda}_{kl} + \boldsymbol{\Lambda}_{ij} ]t)}\int_0^1{\boldsymbol{u'}_{ij}}^{\text{H}}\boldsymbol{u'}_{kl} \diff y,
\end{multline}
where the superscript, `H' denotes Hermitian transpose. Using the orthogonality relations given by Eq.~(\ref{eq:orth}) and normalizations of $\psi_1$ and $\psi_2$, we get,
\begin{equation}
E(t) = \frac{1}{2}  \left ( \sum_{k}|a_k^{}|^2\exp{(2\lambda_{k}t)} + \sum_{l}|a_l^{}|^2 \exp{(2\lambda_{l}t)} \right ) . \label{eq:energy_orth}
\end{equation}
As said in the subsection~\ref{linsys_summary}, in this present case of $\{\alpha = 0, n=0\}$, 
the evolution of perturbations are determined by a dissipative phenomena. Therefore, 
$\lambda_k$ and $\lambda_l$ are all negative. This shows that the energy, $E(t)$ in 
Eq.~(\ref{eq:energy_orth}), is a monotonously decaying function of time. Since the interaction with the mean 
flow is null for these modes, they should be supplied with energy by other means of receptivity 
such as vibration of the wall or by inflow-borne disturbances. The role of acoustic disturbances may 
also be ruled out since extremely longwaves imply extremely low (near-zero) infrasonic frequencies of the source (since, frequency $\sim 1$/wavelength).

\section{Numerical Results
\label{sec:numerical}}

\subsection{Numerical method}
The system of Eqs.~(\ref{eq:phi_in_y})-(\ref{eq:omega_in_y_n_0}) can be written as,
\begin{equation}
Re^{-1} \boldsymbol{A} \boldsymbol{q} = \omega \boldsymbol{B} \boldsymbol{q}, \label{eq:geneigsys}
\end{equation}
where $\boldsymbol{q} = \{ \phi, \Omega \}^T$ and, $\boldsymbol{A}$ and $\boldsymbol{B}$ are matrices of 
operators. The elements of  $\boldsymbol{A}$ and $\boldsymbol{B}$ can be figured out from those 
equations (Eq.~(\ref{eq:phi_in_y})-(\ref{eq:omega_in_y_n_0})) for each case of $n\neq0$ and $n=0$. For 
discretization, we use spectral collocation method in transform space (see, for example, 
Ref.~\cite{trefethen2000spectral} or Appendix A.5 of Ref.~\cite{schmid2001stability} and S. C. Reddy's codes therein). The generalized eigensystem, Eq.~(\ref{eq:geneigsys}) is discretized 
at the points,
\begin{subequations}
\begin{equation}
y_j = 2^{-1}(1+\xi_j), \quad \text{if $\alpha\leq 3$ and $Re \leq 6000$,   } \label{eq:no_stretch}
\end{equation}
else
\begin{equation}
y_j = \left ( \exp(a)-1\right )^{-1}\left ( \exp[a(1+\xi_j)/2]-1\right ), \label{eq:stretch}
\end{equation}
\end{subequations}
where the Gauss Lobatto points, $\xi_j = \cos(\pi j/N), \ j=0\cdots N$ and $N$ is the highest 
order of Chebyshev polynomials. Eq.~(\ref{eq:stretch}) represents a stretching function that 
maps $\xi$ on $y$ such that the grid points cluster towards $y=0$ for the parameter $a>0$. For 
$n  \leq  5$, $a=2$ gives best result, whereas $a=3$ works well for $n>5$.

The stretching is needed in our present formulation for large values of $\alpha$ and $Re$, 
due to the following. As the algebraic $r-$variations of the velocity fields in a power-law 
fashion (i.e., $y-$variations) has been factored out from the unknowns by the Priymak and 
Miyazaki ansatz given in Eq.~(\ref{eq:uvwprimes_r_variations}), the functions
$\phi$ and $\Omega$ do not have any other constraints apart from those of satisfying regularity. 
This allows $\Omega$ and $\phi$ to have steep variations at the pipe centre, given that there is 
a regular singularity at that location. 
The freedom of $\Omega$ and $\phi$ and their derivatives from the requirement to vanish at the 
pipe centre, coupled with their vanishing at the wall and the presence of regular singularity there, 
causes a boundary layer behaviour at pipe centre. 
Hence, the stretching introduced in Eq.~(\ref{eq:stretch}) is required in our formalism.
However, it is precisely this boundary layer behaviour at the pipe centre that allows for 
avoiding the pseudospectra, as will be shown further down.

Note that the requirement for stretching does not arise in the formulation of Burridge \& Drazin, 
where the functions and one of their derivatives vanish at both boundaries, namely, wall and pipe 
centre. Therefore, upon solving their system, the variations of the solutions are spread over the 
entire domain, avoiding steep variations. Such clustering is neither needed in the formulation of 
Meseguer \& Trefethen since the computation is performed with $r$ as independent variable, which 
would allow the collocation to resolve the near-centre of the pipe. 

The no-slip boundary conditions are implemented through a set of row and column operations of reducing 
the order of the system as described through commented lines in the code provided in~\ref{app:code}. We would like to note in passing that the alternate method of spurious modes 
technique, which is a concise and clever way of implementing such homogeneous boundary conditions 
described in Ref.~\cite{schmid2001stability}, renders our results different at the order of $10^{-9}$. 
The method of elimination that we used gives $2N-1$ number of modes in the cases of $n\neq0$ and $\{n=0, \alpha\neq0\}$, and $2N$ modes in the case of $\{n=0, \alpha=0\}$. 
Care is taken to ward off 
round-off errors by evaluating the functions $g_i$, $i=1\cdots7$, which are polynomials in $\alpha^2y$, 
in a nested manner.

The eigensystems were solved using QZ algorithm implemented in the tool, "EIG" of Matlab (version 
R2016a Update 7), and the least decaying modes are evaluated through the Arnoldi method implemented in 
the tool, "EIGS" of Matlab which allows specification of the error tolerance.

\subsection{Validations and comparisons \label{subsec:validations}}

The obtained eigenvalues are compared with that of Meseguer \& Trefethen~\cite{meseguer2000spectral}, 
Schmid \& Henningson~\cite{schmid2001stability} and Priymak \& Miyazaki~\cite{priymak1998accurate} with 
identical parameters as in those references, and shown in Table~\ref{tb:compare_meseg}. It can be noted 
that, on the average, the accuracy matches with that of  Meseguer \& 
Trefethen~\cite{meseguer2000spectral}. 
\begin{table*}\centering
\caption{Least decaying eigenvalues of converged accuracy. `Size' refers to the order of eigensystem\label{tb:compare_meseg}}
{\footnotesize
\begin{tabular}{|lll|l|ll|l|ll|}
\hline
\multicolumn{3}{|c|}{Parameters} & \multicolumn{3}{c|}{Present Method} & \multicolumn{3}{c|}{Meseguer \& Trefethen(2000)}\\
\hline
\multicolumn{1}{|c}{$\alpha$}& \multicolumn{1}{c}{$n$} &\multicolumn{1}{c|}{$Re$}& \multicolumn{1}{l|}{Size}& \multicolumn{1}{c}{$\omega_r$} & \multicolumn{1}{c|}{$\omega_i$} & \multicolumn{1}{l|}{Size}& \multicolumn{1}{c}{$\omega_r$} & \multicolumn{1}{c|}{$\omega_i$} \\
\hline
$1$ & $0$ & $3000$ & $93$ & $0.94836022205056$  & $-0.051973111282766$ & $110$  & $0.9483602220505$ & $-0.0519731112828$\\
$1$ & $1$ & $3000$ & $93$ & $0.9114655676232$  & $-0.041275644694$ & $110$  & $0.91146556762$ & $-0.041275644693$\\
$1$ & $2$ & $3000$ & $93$ & $0.88829765875$  & $-0.060285689555$ & $110$  & $0.88829765875$ & $-0.060285689559$\\
$1$ & $3$ & $3000$ & $93$ & $0.86436392106$  & $-0.083253976943$ & $110$  & $0.86436392104$ & $-0.08325397694$\\
$0$ & $0$ & $3000$ & $94$ & \multicolumn{1}{c}{$0$}  & $-0.001927728654315596$ & $110$  & \multicolumn{1}{c}{$0$} & $-0.0019277286542$\\
$0$ & $1$ & $3000$ & $93$ & \multicolumn{1}{c}{$0$}  & $-0.0048939902144$ & $110$  & \multicolumn{1}{c}{$0$} & $-0.00489399021$\\
$0$ & $2$ & $3000$ & $93$ & \multicolumn{1}{c}{$0$}  & $-0.00879153881$ & $110$  & \multicolumn{1}{c}{$0$} & $-0.0087915387$\\
$0$ & $3$ & $3000$ & $93$ & \multicolumn{1}{c}{$0$}  & $-0.0135688219$ & $110$  & \multicolumn{1}{c}{$0$} & $-0.0135688219$\\
$1$ & $1$ & $9600$ & $99$ & $0.9504813966688$  & $-0.023170795763$ & $110$  & $0.950481396670$ & $-0.023170795764$\\
\hline
\multicolumn{3}{|c|}{Parameters} & \multicolumn{3}{c|}{Present Method} & \multicolumn{3}{c|}{Schmid \& Henningson(2001)}\\
\hline
\multicolumn{1}{|c}{$\alpha$}& \multicolumn{1}{c}{$n$} &\multicolumn{1}{c|}{$Re$}& \multicolumn{1}{l|}{Size}& \multicolumn{1}{c}{$\omega_r$} & \multicolumn{1}{c|}{$\omega_i$} & \multicolumn{1}{l|}{Size}& \multicolumn{1}{c}{$\omega_r$} & \multicolumn{1}{c|}{$\omega_i$} \\
\hline
$1$ & $0$ & $2000$ & $85$ & $0.93675536015933$  & $-0.063745512531531$ & $-$  & $0.93675536$ & $-0.06374551$\\
$0.5$ & $1$ & $2000$ & $77$ & $0.423234848559$  & $-0.0358816618407$ & $-$  & $0.423234850$ & $-0.03588166$\\
$0.25$ & $2$ & $2000$ & $77$ & $0.18137922101$  & $-0.037238251507$ & $-$  & $0.18137922$ & $-0.0372382525$\\
$0$ & $1$ & $2000$ & $73$ & \multicolumn{1}{c}{$0$}  & $-0.0073409853206$ & $-$  & \multicolumn{1}{c}{$0$} & $-0.00734099$\\
\hline
\multicolumn{3}{|c|}{Parameters} & \multicolumn{3}{c|}{Present Method} & \multicolumn{3}{c|}{Priymak \& Miyazaki(1998)}\\
\hline
\multicolumn{1}{|c}{$\alpha$}& \multicolumn{1}{c}{$n$} &\multicolumn{1}{c|}{$Re$}& \multicolumn{1}{l|}{Size}& \multicolumn{1}{c}{$\omega_r$} & \multicolumn{1}{c|}{$\omega_i$} & \multicolumn{1}{l|}{Size}& \multicolumn{1}{c}{$\omega_r$} & \multicolumn{1}{c|}{$\omega_i$} \\
\hline
$20$ & $20$ & $4000$ & $201$ & $1.476280140$  & $-1.0395781217$ & $202$  & \multicolumn{1}{l}{$1.476280140\cdots$} & \multicolumn{1}{l|}{$-1.039578121\cdots$}\\
& & & & & & & \multicolumn{1}{r}{$\cdots 6380943001$} & \multicolumn{1}{r|}{$\cdots 8520833192$} \\
\hline
\end{tabular}
}
\end{table*}
We also found that all 41 eigenvalues listed in~\cite{meseguer2000spectral} were matching at similar 
accuracy with the present computation. The figures for Schmid \& Henningson reported in the 
Table~\ref{tb:compare_meseg} are after converting the complex phase-speeds reported in~\cite{schmid2001stability} into complex $\omega$'s by the relation $\omega = \alpha c$. The 
improved accuracy in the present case is pronounced in comparison with Schmid \& Henningson, which is 
due to the deployment of the ansatz of Priymak \& Miyazaki~\cite{priymak1998accurate} given in 
Eq.~(\ref{eq:vetaprimes_r_variations}). Schmid \& Henningson's~\cite{schmid2001stability} results were 
obtained by solving a system equivalent to that of Burridge \& Drazin~\cite{burridge1969comments}.

As shown in Table~\ref{tb:compare_meseg}, a matching accuracy is reached at collocation points as low as $N +1 = 48$, which corresponds to the system size of $2N-1 = 93$. Increasing the value of $N$ further did not increase the accuracy of the least decaying mode in the present 64-bit calculations, although the 
modes with higher decay rates improved in convergence. 

The attaining of the maximum possible accuracy of the least decaying mode for the case of 
$\alpha = n = 0$ shown in Table~\ref{tb:compare_meseg}, is due to that the system 
is normal, hence the effectiveness of eigenvalue iterative algorithms is maximized. 
In this case, the Arnoldi method reduces to the Lanczos method~\cite{golub2012matrix}.

The eigenvalue of the least decaying mode of the case with parameters, $\alpha = n = 20$ and $Re=4000$ 
matches well with that of Priymak \& Miyazaki~\cite{priymak1998accurate}. It should be noted that the 
present computations have been performed in 64-bit calculations. From the number of significant places 
in the data of Priymak \& Miyazaki shown in Table~\ref{tb:compare_meseg}, one can note that it is a 
result from 80-bit computation.

\begin{figure}\centering
	\includegraphics[width=0.45\textwidth]{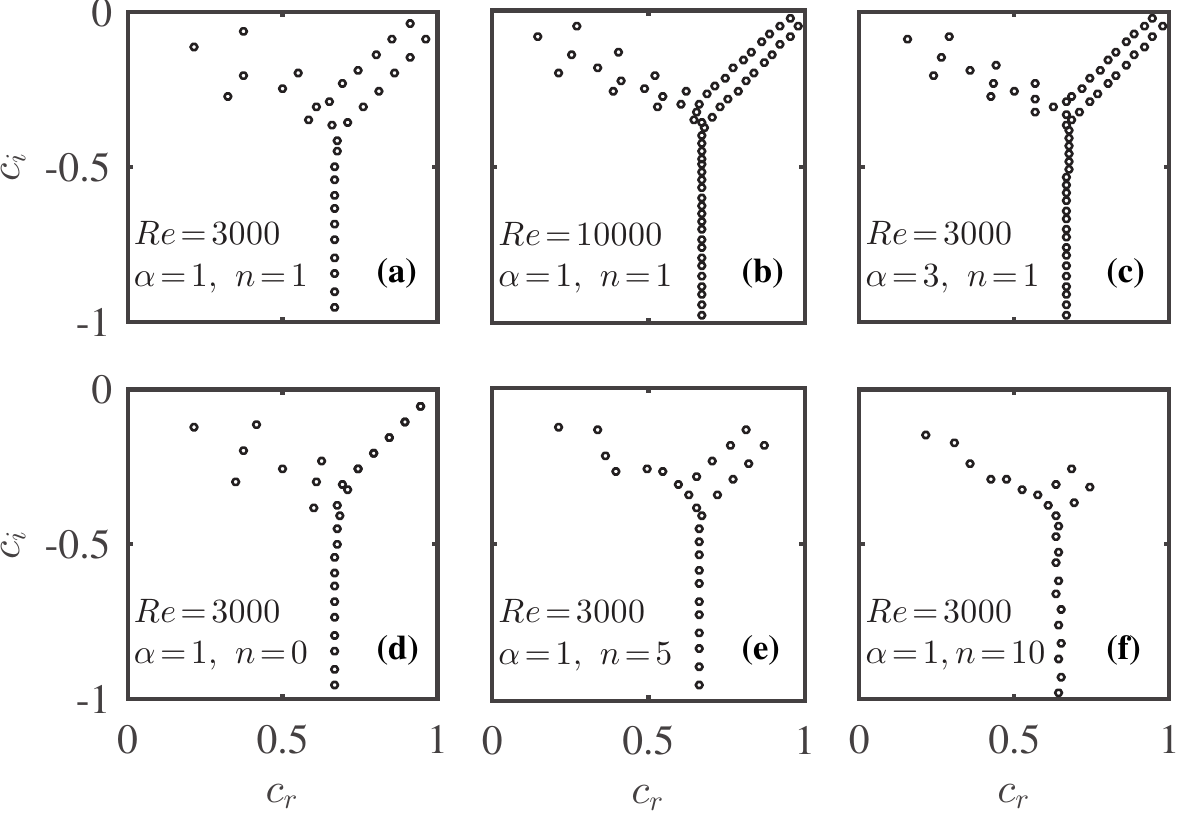}
	\caption{Spectra for parameters shown in each panels. Other parameters:(a) $N=47$; (b) $N=68$; (c) $N=69$; (d) $N=53$; (e) $N=53$; (f) $N=53$;\label{fig:spectra}}
\end{figure}

Figure~\ref{fig:spectra} shows representative spectra for various combinations of parameters. 
The sub-figures show the overall converged spectrum. As can be observed from the caption of 
Fig.~\ref{fig:spectra}, an increase in the number of collocation points is needed when there is 
an increase in $\alpha$, shown in Fig.~\ref{fig:spectra}(c), or $Re$, Fig.~\ref{fig:spectra}(b), 
but seldom for an increase in $n$, shown in Fig.~\ref{fig:spectra}(e) and (f). However, it should 
be observed that there is a surfacing of mild distortion in the spectrum for $n=10$ with 
$\alpha = 1$ and $Re = 3000$ as shown in Fig.~\ref{fig:spectra}(f). This can be explained using Eq.~(\ref{eq:curl:3}). It is well-known that the non-normality is caused by the convective terms that interact with mean shear. These terms are enhanced for large values of $n$ by the first term on the RHS of Eq.~(\ref{eq:curl:3}) giving rise to such distortions.
The spectrum in 
Fig.~\ref{fig:spectra}(b) has same parameters as a sub-figure of Fig.~1 of Meseguer \& 
Trefethen~\cite{meseguer2003linearized} and can be compared as a validation. For the case of 
$\{n=0,\alpha = 1\}$ shown in Fig.~\ref{fig:spectra}(d), the two wall-mode branches appears to 
have merged at the accuracy of the figure. As the azimuthal wavenumber is increased the number 
of modes in the centre branch decrease as shown in Fig.~\ref{fig:spectra}(e) and~(f).

Table~\ref{tb:compare_n_0_al_0} compares the spectrum obtained for the case of 
$\{\alpha =0, n=0\}$ from the system of Eqs.~(\ref{eq:psi1_in_y})-(\ref{eq:psi2_in_y}) with those 
obtained from the characteristic relations given by Eqs.~(\ref{eq:psi1_poly})-(\ref{eq:psi2_poly}). 
These characteristic relations have been approximated by polynomials, as explained in the caption of 
Table~\ref{tb:compare_n_0_al_0}, by setting a cut-off value of 90 for the running index $k$ in the 
summation. 
\begin{table}[htb]\centering
\caption{Comparison of first seven least decaying eigenvalues, $\omega_i$'s for the case of $\{\alpha = 0, n=0\}$ with $Re = 3000$ and $N = 47$ obtained from system of 
Eq.~(\ref{eq:psi1_in_y})-(\ref{eq:psi2_in_y}) versus those obtained from $\omega_i = \lambda Re/4$ 
where the $\lambda$'s are the real solutions of the characteristic 
Eq.~(\ref{eq:psi1_poly})-(\ref{eq:psi2_poly}) approximated with the setting $\max\{k\} = 90$. The 
reported figures are of the converged accuracy with respect to $N$ and $\max\{k\}$.  
\label{tb:compare_n_0_al_0}}
{\footnotesize
\begin{tabular}{|d{20}d{18}|}
\hline \hline \noalign{\smallskip}
\multicolumn{1}{|l}{ $\omega_i$'s from Eq.~(\ref{eq:psi1_in_y})} & \multicolumn{1}{r|}{ $\lambda Re/4$ with $\lambda$'s from Eq.~(\ref{eq:psi1_poly})}\\
\hline
-0.001927728654315596&-0.0019277286543156\\
-0.01015708744788736&-0.010157087447887\\
-0.02496233559689843&-0.024962335596\\
-0.04634676147548659&-0.046346761475\\
-0.0743107678725448&-0.07431076787\\
-0.108854450977443&-0.108854450\\
-0.149977842839345&-0.14997784\\
\hline  \noalign{\smallskip}
\multicolumn{1}{|l}{ $\omega_i$'s from Eq.~(\ref{eq:psi2_in_y})} & \multicolumn{1}{r|}{ $\lambda Re/4$ with $\lambda$'s from Eq.~(\ref{eq:psi2_poly})}\\
\hline
-0.004893990214041297&-0.0048939902140412\\
-0.01640615210723253&-0.01640615210723\\
-0.03449981796504557&-0.0344998179650\\
-0.0591735889379348&-0.05917358893\\
-0.090427218090958&-0.0904272181\\
-0.128260635034236&-0.12826063\\
-0.172673813670569&-0.1726738\\
\hline \hline
\end{tabular}
}
\end{table}
In the table, we have shown only the converged decimal figures that does not change 
when $N$ is changed from $46$ to $47$, or when $\max\{k\}$ is changed from 89 to 90. 
As can be noted from this table, the results obtained by computation from Eqs.~(\ref{eq:psi1_in_y})-(\ref{eq:psi2_in_y}) 
have higher converging precision in comparison with those given by the characteristic relations
Eqs.~(\ref{eq:psi1_poly})-(\ref{eq:psi2_poly}). The characteristic relations contains a quadratic 
term of factorials of $k$ in the denominator, which causes loss of precision for large values of 
$k$. The set value of $\max\{k\} = 90$ is the largest that we could afford in the present 
64-bit (double precision) calculations. Especially, the highly decaying modes are prone to the 
precision loss since $|\lambda|$ are large for these modes. 

\begin{figure}\centering
\includegraphics[width=0.45\textwidth]{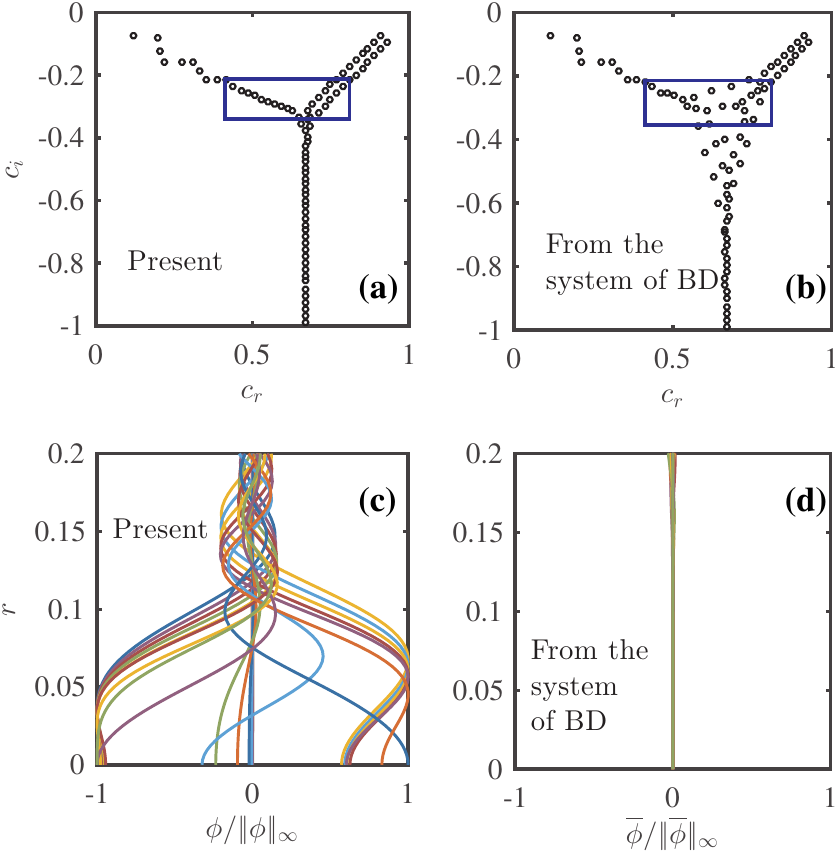}
	\caption{Comparison of spectra and eigenfunctions from the present formulation with that computed from the system of Burridge \& Drazin, Eq.~(\ref{eq:bd_1})-(\ref{eq:bd_2}). Parameters: $N=200$, $\alpha = 10$, $n = 7$, $Re = 2000$. (a) and (b) are spectra; (c) and (d) are the $\infty-$norm normalized real parts of $\phi$ and $\overline{\phi}$, respectively, close to the centreline of 25 modes selected from the region shown by rectangles in (a) and (b). In (c) $\phi$ has been plotted with respect to $r$ for the purpose of comparison with (d). \label{fig:spectra_bd}}
\end{figure}

Finally, in Fig.~\ref{fig:spectra_bd} we compare the result from the present system for large wavenumbers 
with that from solving the system of Burridge \& Drazin~\cite{burridge1969comments}. The system of Burridge 
and Drazin is given by
\begin{align}
(\omega-\alpha U) T\overline{\phi} =& -\alpha d r^{-1} D\left (rd^{-1}U_r\right )\overline{\phi} \nonumber \\
&+ \mi Re^{-1} \left ( T^2 \overline{\phi} -  2\alpha n T \overline{\Omega}\right ) \quad \text{and} \label{eq:bd_1}\\
(\omega-\alpha U) \overline{\Omega} =& -n (rd)^{-1}U_r \overline{\phi} \nonumber \\
&+ \mi Re^{-1} \left ( 2\alpha n d^{-2} T \overline{\phi} +  S \overline{\Omega}\right ), \label{eq:bd_2}
\end{align}
where $D=\frac{\diff}{\diff r}$, $\overline{\phi} = -\mi r v'$ and $\overline{\Omega} = (\alpha r w'-nu')/d$, and the operators $T$ and $S$ are defined as $T[\cdot] = dr^{-1}D(rd^{-1}D[\cdot]) - dr^{-2}[\cdot]$ and $S[\cdot] = (rd)^{-1}D(rdD[\cdot]) - dr^{-2}[\cdot]$. The boundary conditions are 
given in Schmid \& Henningson~\cite{schmid2001stability}. (However, there is a trivial typographical 
sign error in the first term of Eq.~(3.41) in that reference. The system of Burridge \& Drazin as it 
appear in Ref.~\cite{burridge1969comments}, has a typographical error in the definition of operator, 
$T$.) 

As can be noted from Fig.~\ref{fig:spectra_bd}(a) and (b), the computation from the system, 
Eq.~(\ref{eq:bd_1})-(\ref{eq:bd_2}), suffers from pseudospectra, which is absent in the present 
formulation. This can be understood as follows. In Fig.~\ref{fig:spectra_bd}(a) and (b) the blue 
rectangle shows a region of the spectrum containing 25 modes that undergoes distortion when 
computation is performed from the system of Burridge and Drazin. Fig.~\ref{fig:spectra_bd}(c) 
and (d) shows $\phi$'s of the modes that are present within those rectangles.  As mentioned earlier, the 
components of eigenfunctions of the present formulation, $\phi$ for example, undergoes steep 
variations close to the centreline for large $\alpha$. Such variations have been suppressed in 
the case of the system of Eq.~(\ref{eq:bd_1})-(\ref{eq:bd_2}), as 
$\overline{\phi} = -\mi r^n \phi$, which shows that Burridge and Drazin's unknown 
$\overline{\phi}$ is a multiplication of the present unknown $\phi$ by a function that is 
vanishingly small for a range of the domain. When $n$ is large, $\overline{\phi}$ is closer 
to zero for almost a fifth of the pipe radius as can be seen from Fig.~\ref{fig:spectra_bd}(d) for the modes of distorted regions of the spectrum shown in
Fig.~\ref{fig:spectra_bd}(b). This implies that the details of the factor of $\phi$ in the expression for $\overline{\phi} (=-\mi r^n \phi)$, 
where the former, i.e., $\phi$ is the distinguishing feature among different eigenfunctions, 
would be poorly represented even at double precision. For evidence that the distortion of the spectrum in 
Fig.~\ref{fig:spectra_bd}(b) originiates in round-off errors, one can note that this pattern is similar to Fig.~2(b) of Ref.~\cite{schmid1994optimal}, which shows a spectrum in 16 bit, 
although for a set of small values of parameters, i.e., $\alpha$, $n$ and $Re$.

Let us assume that we cast Eqs.~(\ref{eq:bd_1})-(\ref{eq:bd_1}) in the form of 
$Re^{-1} \boldsymbol{\overline{A}} \boldsymbol{\overline{q}} = \omega \boldsymbol{\overline{B}} \boldsymbol{\overline{q}}$ where $\overline{\overline{q}} = (\overline{\phi}, \overline{\Omega})^{\rm T}$, similar to the present formulation as we saw in Eq.~(\ref{eq:geneigsys}). Let us compare the 
condition numbers, $\text{cond}(\boldsymbol{L}-\omega I)$ and $\text{cond}(\boldsymbol{\overline{L}}-\omega I)$ where $\boldsymbol{L} = \boldsymbol{B}^{-1}\boldsymbol{A}/Re$ and  $\boldsymbol{\overline{L}} = \boldsymbol{\overline{B}}^{-1}\boldsymbol{\overline{A}}/Re$. We found that $\text{cond}(\boldsymbol{\overline{L}}-\omega I)\sim O(10^{18})$ and $\text{cond}(\boldsymbol{L}-\omega I)\sim O(10^{21})$, for $\omega's$ chosen from their respective rectangles shown in 
Fig.~\ref{fig:spectra_bd}(a) and (b). This shows that the present eigensystem is much closer to 
singularity for these modes than that of Burridge and Drazin. In other words, the contours of resolvent 
norms in the complex plane will be of smaller value for the same radial distance from these eigenvalues in 
the present formulation in comparison with the system of Burridge and Drazin.

\section{Inviscid Algebraic Instability of Axially constant Modes \label{sec:ellingsen}}
In this section we deduce the present flow configuration's analogue of the findings of 
Ellingsen and Palm~\cite{ellingsen1975stability} for plane shear flows for the streamwise constant modes. The aim is to find 
the initial $\phi$ and $\Omega$, i.e. at time $t=0$, that maximizes the energy amplification, and the output pattern at a
later time $t$. Ellingsen and Palm solution concerns inviscid solutions with $\alpha =0$ in a non-modal setting and results 
in a initial value problem. (For such initial value problem in viscous scenario, see Ref.~\cite{bergstrom1992initial} for 
series solutions, and Ref.~\cite{schmid1994optimal} for optimal patterns). For axially constant inviscid modes that has a 
general non-modal time dependence, Eq.~(\ref{eq:phi_in_y}) and~(\ref{eq:omega_in_y}) become
\begin{align}
\frac{\partial }{\partial t}\left [ (\ell+2)\mathcal{D}\phi + y\mathcal{D}^2\phi\right ] =& 0 \label{eq:inv_nonmodal_phi}\\
\frac{\partial \Omega}{\partial t} =& -2\mi n U_y \phi \label{eq:inv_nonmodal_omega}.
\end{align}
The Eq.~(\ref{eq:inv_nonmodal_phi}) implies conservation of kinetic energy in the radial and azimuthal 
directions. As will be shown later in this section, the sum of energies in these directions is 
proportional to $(\ell+2)\mathcal{D}\phi + y\mathcal{D}^2\phi$ (apart from a factor of $-2n^{-2}y^{\ell+1}$). 
Eq.~(\ref{eq:inv_nonmodal_phi}) governs the rate at which energy is pumped into the radial vorticity 
by the mean shear. Since the azimuthal velocity is conserved, this implies the rate at which the axial 
component of the kinetic energy is enhanced, giving rise to streaks.

Eq.~(\ref{eq:inv_nonmodal_phi}) can be readily integrated to yield $(\ell+2)\mathcal{D}\phi + y\mathcal{D}^2\phi = C(y)$, 
where the function, $C(y) = [(\ell+2)\mathcal{D}\phi + y\mathcal{D}^2\phi]|_{t=0}$ is equivalent to the specification of 
an initial value of the kinetic energy in the non-streamwise directions.  Such specification of $C(y)$ at $t=0$ is 
equivalent to the specification of $\phi_0$, the initial value of $\phi$ itself. (In fact, for a specified $C(y)$, the 
$\phi_0(y)$ that is finite at centreline and satisfying $\phi_0(1) = 0$ is given by 
$\phi_0(y) = \int_0^1  [\mathcal{U}(y-\overline{y}) -1 ] {\overline{y}}^{-(\ell+2)}\int_0^{\overline{y}} C(\tilde{y}) \tilde{y}^{\ell+1} \diff \tilde{y} \diff \overline{y}$ 
where $\mathcal{U}(y)$ is the unit step function). Therefore,
\begin{align}
\phi(t,y) =& \phi_0(y) \label{eq:inv_nonmodal_phi_soln}\\
\Omega(t,y) =& \Omega_0(y) -2\mi n U_y t \phi_0(y), \label{eq:inv_nonmodal_omega_soln}
\end{align}
where $\Omega_0(y) = \Omega(t=0,y)$. Eq.~(\ref{eq:inv_nonmodal_omega_soln}) is obtained from integration of 
Eq.~(\ref{eq:inv_nonmodal_omega}) and
describes algebraic growth since the second term is secular, which implies the collapse of linear perturbation theory after certain initial time unless viscosity acts and kills it. In terms of velocities these equations translate into: $v'(t,y) = v'(0,y)$, $u'(t,y) = u'(0,y) - 2t\sqrt{y}U_yv'(0,y)$ and $w'(t,y) = \mi n^{-1} [v'(0,y) + 2y\mathcal{D}v'(0,y)]$.

To find the optimal patterns in this inviscid limit, we follow the method of Ref.~\cite{hanifi1998compressible} 
for compressible flow with necessary modifications for the present incompressible configurations 
(see also Ref~\cite{malik2008linear}). As we are working with a 2-tuple variable, the constraint of continuity is 
already taken into consideration. The energy, $E(t) \equiv2^{-1}\int_0^1 |\boldsymbol{u'}|^2dy$ can be written with 
the help of Eq.~(\ref{eq:vetaprimes_r_variations})-(\ref{eq:wprime_new}) as 
\begin{multline}
E(t) = \frac{1}{2n^2}\int_0^1 \left [ y^{\ell+1}|\Omega|^2 + 2n^2y^\ell|\phi|^2 + 4y^{\ell+2}|\mathcal{D}\phi|^2  \right . \\ \left . + 2(\ell+1)y^{\ell+1} \left ( \phi^*\mathcal{D}\phi + \phi\mathcal{D}\phi^* \right ) \right ] \diff y,
\end{multline}
which can be written after integration by parts as, 
\begin{equation}
E(t) = (2n^2)^{-1}\int_0^yy^{\ell+1} \boldsymbol{q_0^\text{H}\mathcal{A}^\text{H}\mathcal{MA}q_0} \diff y,
\end{equation}
where $\boldsymbol{q_0} = \{\phi_0,\Omega_0\}^{\textrm{T}}$, the positive definite 
$2\times 2$ matrix $\boldsymbol{\mathcal{M}} = \textrm{diag}\{-4[(\ell+2)\mathcal{D}+y\mathcal{D}^2],1\}$. 
The operator $\boldsymbol{\mathcal{A}}$ is the propagator of $\boldsymbol{q_0}$ and can be obtained by casting 
Eq.~(\ref{eq:inv_nonmodal_phi_soln}) and Eq.~(\ref{eq:inv_nonmodal_omega_soln}) in matrix form. The elements of 
$\boldsymbol{\mathcal{A}}$ are $\mathcal{A}_{11} = \mathcal{A}_{22} = 1$, $\mathcal{A}_{12} = 0$, and 
$\mathcal{A}_{21} = -2\mi n U_yt$. Let $g(t)$ be the propagator of the initial energy, i.e., 
\begin{equation}
E(t) = g(t)E(0). \label{eq:ege}
\end{equation} 
Note that $E(0) = (2n^2)^{-1}\int_0^yy^{\ell+1} \boldsymbol{q_0^\text{H}\mathcal{M}q_0} \diff y$. Let us define, $G(t) = \max_{\boldsymbol{q_0}}g(t)$. This $G(t)$ is found as follows. Taking the functional derivative of Eq.~(\ref{eq:ege}), we get 
\begin{equation}
\frac{\delta }{\delta \boldsymbol{q_0}^\text{H}} E(t) = \left (\frac{\delta }{\delta \boldsymbol{q_0}^\text{H}}g(t)\right ) E(0) + g(t) \frac{\delta }{\delta \boldsymbol{q_0}^\text{H}} E(0).
\end{equation}
Setting $\delta g/\delta \boldsymbol{q_0}^\text{H} = 0$ for maximization we get,
\begin{equation}
\boldsymbol{\mathcal{A}^\text{H}\mathcal{MA}q_0} = g(t)\boldsymbol{\mathcal{M}q_0}, \label{eq:inviscid_g_eig_eqn}
\end{equation}
which is a differential equation that can also be written as
\begin{multline}
2 \left [ (nU_yt)^2-(\ell+2)\mathcal{D} -y\mathcal{D}^2 \right ]\phi_0 + \mi nU_yt \Omega_0 \\ = -2g(t)\left [ (\ell+2)\mathcal{D} +y\mathcal{D}^2 \right ] \phi_0
\end{multline}
\begin{equation}
-2\mi nU_yt\phi_0 + \Omega_0 = g(t) \Omega_0 \label{eq:omg_opt}
\end{equation}
\begin{figure}[htb]
\includegraphics[width=0.45\textwidth]{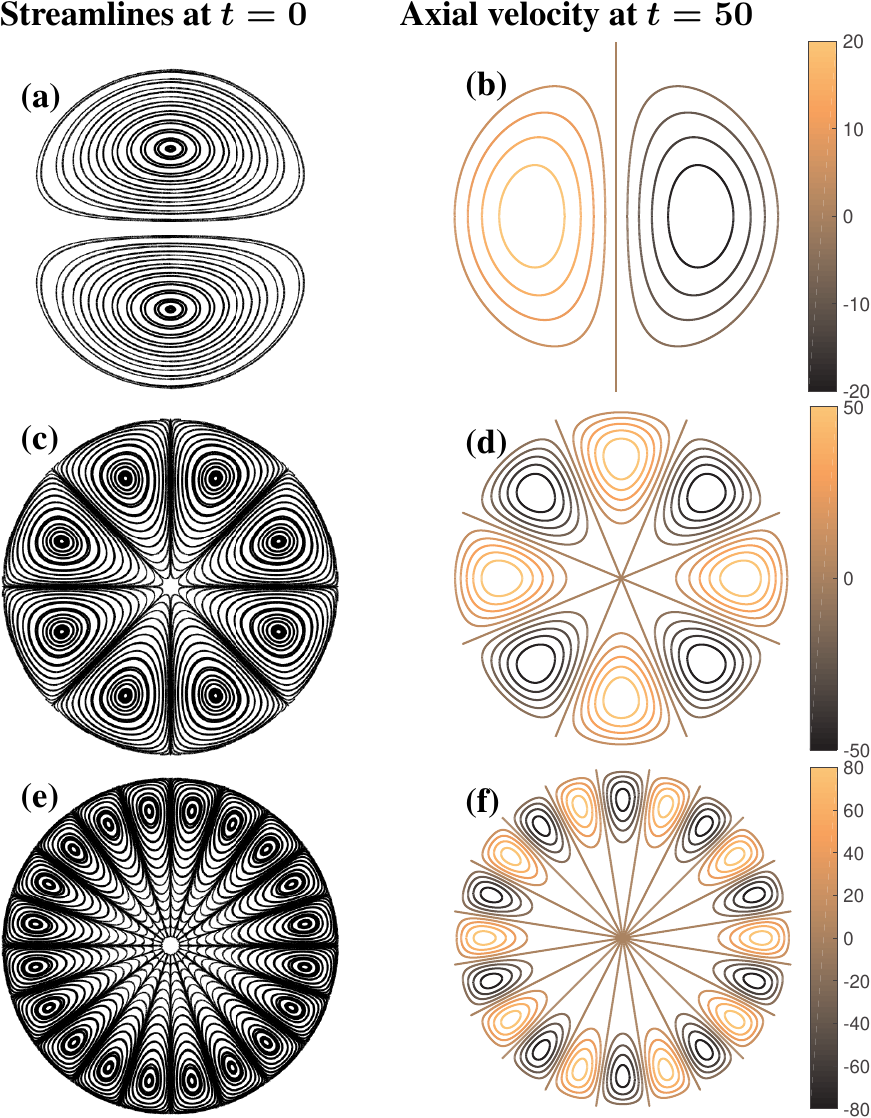}
	\caption{The inviscid initial (left) and the instantaneous velocity patterns at a later time (right) for axially constant modes that are optimized for maximum energy growth: (a) and (b) $n=1$; (c) and (d) $n=4$; and, (e) and (f) $n=10$. The sub-figures on the left column are the 2D streamlines in the cross section of the pipe at time, $t=0$. The sub-figures on the right column are the contours of the axial velocity at time, $t = 50$ with levels indicated by colour-bar. \label{fig:opt_invisc}}
\end{figure}
This is an eigenvalue problem with the Lagrange multiplier, $g(t)$, as the eigenvalue for a chosen $t$. Since the operators 
are Hermitian, $g(t)$'s are all real, and as $\mathcal{M}$ is positive definite by definition, $g(t)$'s are all positive. 
The boundary conditions are no-penetration at the wall and regularity at the centreline, i.e., $\phi_0(1) = 0$ and 
$2[(g(t)-1)(\ell+2)\mathcal{D} +(nU_yt)^2 ]_{y=0}\phi_0 + \mi nU_yt\Omega_0(0) = 0$, respectively. As $\Omega_0$ appears algebraically in Eq.~(\ref{eq:omg_opt}), it does not require boundary conditions. In fact, we can explicitly solve for $\Omega_0$ in terms 
of $\phi_0$ and express Eq.~(\ref{eq:inviscid_g_eig_eqn}) as an eigenvalue problem with only $\phi_0$ as the unknown. 
However, this will lead to that the eigenvalue parameter, $g(t)$, appears quadratically and will warrant a companion matrix 
method to solve the equation numerically, which is equivalent to the problem as it appears currently in 
Eq.~(\ref{eq:inviscid_g_eig_eqn}). The maximum amplification, $G(t)$ as defined above, is given by $G(t) = \max_j g_j(t)$ 
where $g_j(t)$'s for a chosen $t$ are the set of eigenvalues of Eq.~(\ref{eq:inviscid_g_eig_eqn}). 
However, the variation $G(t)\sim t^2$ can be anticipated from Eq.~(\ref{eq:inv_nonmodal_omega_soln}). 
In addition, such a variation of $G(t)$ with respect to $t$ has been observed in wall-bounded plane shear 
flows~\cite{gustavsson1991energy}. 
Here, we focus is on the optimal $\boldsymbol{q_0}$ and the corresponding $\boldsymbol{q}(t)$ as obtained from 
Eq.~(\ref{eq:inviscid_g_eig_eqn}).
To solve the system given by Eq.~(\ref{eq:inviscid_g_eig_eqn}), we adopt the the spectral method described in 
previous section.

Figs.~\ref{fig:opt_invisc}(a), (c) and (e) show the streamlines of initial velocity pattern in the cross sectional plane for 
$n=1, 4 \  \text{and} \ 10$, respectively, that are prone to the maximum energy amplification. These patterns represent 
counter rotating vortices. The explanation for the appearance
of these patterns is similar to the case of plane shear flows, such as boundary layer, Couette and plane Poiseuille flows. 
Since counter rotating vortices contain movement of fluid particles in the radial direction,  they can transfer energy 
from mean shear into axial velocity. This phenomenon is widely known as \emph{lift-up effect}~\cite{schmid2001stability} 
in these plane shear flows. In the context of this pipe flow, the term, \emph{up} refers to the radial direction.

Such patterns have been observed in this flow configuration through DNS of full nonlinear viscous equations 
when the mean flow was perturbed initially with axially constant modes~\cite{zikanov1996instability}. 
Fig.~\ref{fig:opt_invisc}(a) matches with the optimal pattern from the viscous flow computation at $Re = 3000$ shown in 
Ref.~\cite{meseguer2003linearized}. Since this inviscid analysis captures those features, one can conclude that the inviscid 
lift-up effect is dominant at initial times before viscosity break this vortices into small scales. 
Figs.~\ref{fig:opt_invisc}(b), (d) and (f) show the contours of axial component of perturbation velocity at a 
later time chosen as $t = 50$ in the non-dimensional scale, where the initial state $\boldsymbol{q_0}$ has been normalized 
such that $E(0) = 1$. These patterns describe near-wall streaks when superimposed with the mean 
flow. These streaks are known to be the characteristic feature of \emph{bypass transition}. (For a matching of these optimal 
perturbations with experimental data in the case of boundary layers, see~\cite{luchini2000reynolds}). As the value of $n$ is 
increased, these streaks are pushed closer to the wall, a result tied to the following two facts: (1) The streamwise 
velocity's growth is proportional to the radial velocity as can be seen from the relation, 
$u'(t,y) = u'(0,y) - 2t\sqrt{y}U_yv'(0,y)$; and, (2) the radial velocity culminates at the periphery of the cross-section 
when $n$ is increased as can be to noted from $v' = r^l \phi_0$.

It should be noted that the $G(t)$ at a chosen time would increase with an increase in $n$ as evident from 
Eq.~(\ref{eq:inv_nonmodal_omega_soln}). This also shows up on the range of the colour bars in Fig.~\ref{fig:opt_invisc} 
as we move down from panels (b) to (f). However, as can be seen from the optimal patterns (in the inviscid limit), 
the size of these vortices reduces when $n$ is increased, making them more prone to the action of viscosity to break them 
further and eventually being killed. This phenomenon can be observed from Fig.4(a) of Ref.~\cite{schmid1994optimal} and 
Fig.4 of Ref.~\cite{bergstrom1993optimal}.

\section{Conclusion}
The present formulation, which uses the ansatz of Priymak \& Miyazaki~\cite{priymak1998accurate} 
results in accurate spectra even for reasonably large values of the parameters, namely, the 
wavenumbers and Reynolds number. The obtained accuracy is on par with that of Meseguer \& 
Trefethen~\cite{meseguer2000spectral}. However, as the boundary conditions of the wall are 
directly imposed on the unknowns used, which is in contrast to that of Meseguer \& 
Trefethen~\cite{meseguer2000spectral}, where they are imposed on the basis functions, the present 
formulation enables the derived system to be applied on a range of flow configurations with 
various boundary conditions. 

The multiplicity of the poles of regular singularity yielded the necessary number of regularity 
conditions in order to complement the boundary conditions at wall, providing 6 conditions to solve the 
6th order system for the cases other than that of $\{\alpha=0, n=0\}$. 

In the case of  $\{\alpha=0, n=0\}$, the working variables had to be changed from a 
representative of normal vorticity and normal velocity to that of streamwise and azimuthal 
velocities. The obtained self-adjoint system's analytical characteristic relation predicted eigenvalues 
which were used to validate the numerical results for this case. The orthogonality showed that the 
superposition of these modes are always decaying.

As a demonstration of the performance of this formulation, the results were contrasted against those from the system of 
Burridge and Drazin~\cite{burridge1969comments}. The performance of the current formulation is found to be highlighted 
at large values of the parameters, and found to postpone the appearance of pseudo-spectra to the regimes of further 
higher values of these parameters. This was tracked to the fact that the numerical representation of distinctive 
features of different eigenfunctions are enhanced in the present formulation, and that it undergoes precision loss in the 
case of Burridge and Drazin due to round-off errors.

Finally, the inviscid algebraic (nonmodal) growth of the perturbation for streamwise constant modes were studied. 
The result was the pipe flow's counterpart of Ellingsen and Palm solutions for plane shear flows. The present working 
variables for the unknowns and the independent variable, $y$, resulted in the solutions which were almost identical in 
structure to that of plane shear-flows, hence retaining the well known characteristics such as having the counter rotating 
vortices and streaks as optimal patterns. These streaks were found to be closer to the wall when the azimuthal wavenumber is larger. For the appearance of such streaks, it is not necessary that the infinitesimal perturbations acquired through receptivity should have the Fourier components with $n$ very large. The nonlinearity can play the role even at very early stages to give birth to such modes through convolution of modes with lower $n$, and subsequently amplified by the linear mechanism of \emph{lift-up}.

\appendix
\section{Operators and functions  \label{app:oper}}
To facilitate locating the functions and operators below, their names have been given in boldface.
\noindent $\boldsymbol{L_1} = (f_1 + f_2 D + rd^{-1}D^2)$, $\boldsymbol{L_2} = (f_3+d^{-1}D)$, $\boldsymbol{L_3} = (f_4 + f_5D +\alpha^2rUd^{-1}D^2)$, $\boldsymbol{L_4} = d^{-1}(f_6 + UD)$, $\boldsymbol{L_5} = (f_7 + f_8D + f_9D^2 + f_{10}D^3 + rd^{-1}D^4)$, $\boldsymbol{L_6} = (f_{11} + f_{12}D + f_{13}D^2 +d^{-1}D^3)$, $\boldsymbol{L_7} = f_{14} + f_{15}D + (rd)^{-1}D^2$, $\boldsymbol{L_8} = f_{16} + d^{-1}D$, $\boldsymbol{L_9} = f_{17} + U_r(rd)^{-1}D$, $\boldsymbol{L_{10}} = d^{-1}(U_r + UD)$, $\boldsymbol{L_{11}} = f_{18} + f_{19}D + f_{20}D^2 + f_{21}D^3 + (rd)^{-1}D^4$, $\boldsymbol{L_{12}} = f_{22} + f_{23}D + f_{24} D^2 + d^{-1}D^3$, $\boldsymbol{L_{13}} = f_{25} - f_{26}D -d^{-1}D^2$, $\boldsymbol{L_{14}} = f_{27} + f_{28}D + D^2$, $\boldsymbol{L_{15}} = f_{29} + f_{30}D + f_{31}D^2+f_{32}D^3 + r^2D^4$, $\boldsymbol{L_{16}} = f_{33} + f_{34}D + r^2d^{-1}D^2$ \\

\noindent $\boldsymbol{f_1} = [(\ell+1)d_2-d^2]/(d^2r)$, $\boldsymbol{f_2} = [d_2 + (\ell+2)d]/d^2$, $\boldsymbol{f_3} = d_2/(d^2r)$, $\boldsymbol{f_4} = [U_r(\alpha^2 r^2 -n^2\ell)-rdU_{rr}]/(dr^2)$, $\boldsymbol{f_5} = -n^2U_r/(rd)$, $\boldsymbol{f_6} = (Ud_2+rdU_r)/(rd)$, $\boldsymbol{f_7} = [(\ell+1)(d_4 + dd_3-(d+1)d_2d^2)+(d+1-\ell^2)d^4]/(d^4r^3)$, $\boldsymbol{f_8} = [d_4+(3\ell + 7)dd_3+(2\ell+ 3-d)d^2d_2 -3(\ell+1)d^4 -(\ell+2)d^3]/(d^4r^2)$, $\boldsymbol{f_9} = [3d_3 +(3\ell+11)dd_2 + (\ell+2-2d)d^2]/(rd^3)$, $\boldsymbol{f_{10}} = [3d_2+(\ell+5)d]/d^2$, $\boldsymbol{f_{11}} = [d_4 + dd_3-(d+1)d^2d_2+2d^4]/(d^4r^3)$, $\boldsymbol{f_{12}} = [3d_3+2dd_2-(d+1)d^2]/(d^3r^2)$, $\boldsymbol{f_{13}} = [3d_2+d]/(d^2r)$,\\ 

\noindent $\boldsymbol{f_{14}} = [(\ell+1)(d+d_5)-d^2]/(d^2r^3)$, $\boldsymbol{f_{15}} = [d_5+(\ell+3)d]/(d^2r^2)$, $\boldsymbol{f_{16}} = (d_5+3d)/(rd^2)$, $\boldsymbol{f_{17}} = (\ell + 1)U_r/(r^2d)$, $\boldsymbol{f_{18}} = (\ell+1)(r^5d^4)^{-1}[d_7 + 2dd_6 - (d+1)d^2d_5 + (1-d)d^3] + r^{-5}(d - \ell^2 + 2\ell -1)$, $\boldsymbol{f_{19}} = [d_7 + (3\ell + 8)dd_6 + (4\ell +7 -d)d^2d_5 -(\ell+1)d^3 - (3\ell+2)d^4]/(rd)^4$, $\boldsymbol{f_{20}} = [3d_6 + (3\ell+13)dd_5 +(2\ell+5)d^2-2d^3]/(rd)^3$, $\boldsymbol{f_{21}} = [3d_5 +(\ell+6)d]/(rd)^2$, $\boldsymbol{f_{22}} = [d_7 + 8dd_6 + (13-d)d^2d_5 + 3(1-d)d^3]/(d^4r^3)$, $\boldsymbol{f_{23}} = (3d_6 + 15dd_5 +13d^2 - d^3)/(d^3r^2)$, $\boldsymbol{f_{24}} = (3d_5 + 8d)/(d^2r)$,\\

\noindent $\boldsymbol{f_{25}} = [d^2-(\ell+1)d_2]/(rd)^2$, $\boldsymbol{f_{26}} = [d_2+(\ell+2)d]/(d^2r)$, $\boldsymbol{f_{27}} = [(\ell^2-d-1)d^2 + 2n^2d_2]/(rd)^2$, $\boldsymbol{f_{28}} = [(2\ell+1)d+2n^2]/(rd)$, $\boldsymbol{f_{29}} = [(\ell+1)(d_4-d^3d_2-d_1d_3+dd_1d_2)+(d-\ell^2)d^4+(2\ell n^2-d_1)d^3]/(d^3r^2)$, $\boldsymbol{f_{30}} = [d_4 + 2(2\ell+3)dd_3 -d_2d_3 -2(2\ell+1)d^4 +(\ell+2)d^2d_1-(2\ell+3)dd_1d_2]/(d^3r)$, $\boldsymbol{f_{31}} = [3d_3+(5\ell+9)dd_2 -2d_2^2-(\ell+2)dd_1 -2d^3]/d^2$, $\boldsymbol{f_{32}} = 2r[d_2+(\ell+2)d]/d$, $\boldsymbol{f_{33}} = (d_3+dd_2-d^3)/d^3$, $\boldsymbol{f_{34}} = r[2d_2+d]/d^2$\\ 

\noindent $\boldsymbol{d_1} = n^2-\alpha^2r^2$, $\boldsymbol{d_2} = (\ell+1)d-2\alpha^2r^2$, $\boldsymbol{d_3} = (\ell d-4\alpha^2r^2)d_2 + 2(\ell-1)\alpha^2r^2d$, $\boldsymbol{d_4} = (\ell-1)dd_3-6\alpha^2r^2d_3+2\alpha^2r^2d[(\ell-4)d_2+(\ell-1)(\ell+2)d-2(\ell-1)\alpha^2r^2]$, $\boldsymbol{d_5} = \ell d-2\alpha^2r^2$, $\boldsymbol{d_6} = (\ell-1)dd_5 -4\alpha^2r^2d_5 +2(\ell-2)\alpha^2r^2d$, $\boldsymbol{d_7} = (\ell-2)dd_6-6\alpha^2r^2d_6+2\alpha^2r^2d[(\ell-5)d_5+(\ell-2)(\ell+1)d-2(\ell-2)\alpha^2r^2]$\\

\noindent $\boldsymbol{g_1} = n^6 + 2(\ell+1)n^4 + n^2[3n^2+4(\ell+1)]\alpha^2y + [3n^2+2(\ell+1)]\alpha^4y^2 + \alpha^6y^3$, $\boldsymbol{g_2} = 4(\ell+2)n^6 + 4n^4(3\ell+5)\alpha^2y + 4n^2(3\ell+4)\alpha^4y^2 + 4(\ell+1)\alpha^6y^3$, $\boldsymbol{g_3} = n^2(\ell^4 + 8\ell^3 + 26\ell^2 + 32\ell + 13) + (3\ell^4 + 20\ell^3 + 50\ell^2 + 36\ell + 3)\alpha^2y + (3\ell^2 + 10\ell + 7)\alpha^4y^2 + \alpha^6 y^3$, $\boldsymbol{g_4} = n^4(\ell^3 + 6\ell^2 + 11\ell + 6) + n^2(3\ell^3+ 15\ell^2 + 21\ell + 5)\alpha^2y + (3\ell^3 + 12\ell^2 + 13\ell + 4)\alpha^4y^2 + (\ell+1)\alpha^6 y^3$, $\boldsymbol{g_5} = 8[n^6(2\ell^2 + 10\ell + 12) + n^4(5\ell^2 + 22\ell + 21)\alpha^2y + n^2(3\ell^2 + 12\ell + 9)\alpha^4y^2 -n^2\alpha^6 y^3 - \alpha^8y^4]$, $\boldsymbol{g_6} = n^2 +2(\ell+1) + \alpha^2y$, $\boldsymbol{g_7} = n^2(\ell^2 + 4\ell + 7) + 2(\ell^2 + 3\ell + 2)\alpha^2y + \alpha^4y^2$, 
$\boldsymbol{g_8} = n^2[3n^2 +4(\ell+1)]$, $\boldsymbol{g_9} = n^4[2(5\ell+9)-n^2]$, $\boldsymbol{g_{10}} = 4n^6(\ell+3)$, $\boldsymbol{g_{11}} = n^4[n^2 + 2\ell + 18]$, $\boldsymbol{g_{12}} = 3\ell^4 + 20\ell^3 + 50\ell^2 + 36\ell + 3$, $\boldsymbol{g_{13}} = n^2(\ell^4 -16\ell^3 -94\ell^2 -136\ell -27)$, $\boldsymbol{g_{14}} = 8n^4(\ell^3+\ell^2-11\ell-15)$, $\boldsymbol{g_{15}} = 16n^6(\ell^2+7\ell+12)$
\section{MATLAB Code \label{app:code}}
{
\begin{lstlisting}
%%%%%%%%%%%%%%%%%%%%%%%%%%%%%%%%%%%%%%%%%%%%%%%%%%%%%%%%%%%%%%%%%%%%%%%%%%%
%                    LINEAR STABILITY CODE FOR PIPE FLOW                  %
%                                                                         %
%                          M Malik & Martin Skote,                        %
%     A linear system for pipe flow stability analysis allowing for       %
%     boundary condition modifications, Computers and Fluids (2019),      %
%            https://doi.org/10.1016/j.compfluid.2019.104267              %
%%%%%%%%%%%%%%%%%%%%%%%%%%%%%%%%%%%%%%%%%%%%%%%%%%%%%%%%%%%%%%%%%%%%%%%%%%%
%{
N      - Highest order of Chebyshev polynomials used. N=47 is best for
         al=0 to 1 and Re < 3000.
al     - Axial wavenumber alpha as in exp[i(alpha x + n theta) - i omega t]
n      - Azimuthal wavenumber as in exp[i(alpha x + n theta) - i omega t]
Re     - Reynolds number
s_bool - Stretching enabler. If s_bool=1, clusters grids towards y=0, else
         grids from -1<xi<1 is mapped linearly on 0<=y<=1.
         Advisable to enable only for al>3 or Re>6000.
a      - Stretching parameter. a=2 gives best result for al>3 or Re>6000
         for small n. a=3 is best for n>6.
xi     - Gauss Lobatto grid points. (-1<=xi<=1)
y      - r^2 (square of radial coordinate) (0<=y<=1)
D0     - Cheby. poly. matrix (maps from coeff. to Gauss Lobatto(GL) pts)
D1-D4  - Differentiation matrices. Maps Cheb coeffs. on 0<=y_j<=1
U,Uy   - Mean velocity,U(y), and its derivative with respect to y
A,B    - Operators in the equation: A q = e Re B q
e    - Complx frequencies omega's as in exp[i(alpha x + n theta - omega t)]
q      - Eigenfunctions.  To note: phi(y) or psi_1(y) = D0*q(1:N+1,:);
         Omega(y) or psi_2(y) = D0*q(N+2:2N+2,:) (For definition of phi,
         Omega, psi_1 and psi_2, see the paper by the authors)
%}
close all; clear variables
%--------------------------------------------------------------------------
% Control variables are N, al, n, Re, s_bool and a only.
N = 200; al = 10; n =7; Re = 2000; s_bool = 1; a=3.0;
%--------------------------------------------------------------------------
l=abs(n)-1; N1=N+1; j=(0:N)'; xi=cos(pi*j/N); D1xi=zeros(N1); D2xi = D1xi;
D3xi = D1xi; D4xi = D1xi; ons=ones(1,N1); D0=cos((j*(pi*j'))/N);
D1xi(:,2:3) = [D0(:,1) 4*D0(:,2)]; D2xi(:,3) = 4*D0(:,1);
for ord=3:N%D^kT_n=2nD^{k-1}T_{n-1}+[n/(n-2)]D^kT_{n-2}, D^k
   jj = ord+1; D1xi(:,jj)=2*ord*D0(:,jj-1)+ord*D1xi(:,jj-2)/(ord-2);
   D2xi(:,jj)=2*ord*D1xi(:,jj-1)+ord*D2xi(:,jj-2)/(ord-2);
   D3xi(:,jj)=2*ord*D2xi(:,jj-1)+ord*D3xi(:,jj-2)/(ord-2);
   D4xi(:,jj)=2*ord*D3xi(:,jj-1)+ord*D4xi(:,jj-2)/(ord-2);
end;
if s_bool==1, y=(exp(a*(xi+1)/2)-1)/(exp(a)-1);
  xiyf = (exp(a)-1)./(1+(exp(a)-1)*y); xi_y=(2/a)*xiyf;
  xi_yy=-(2/a)*xiyf.^2; xi_yyy=(4/a)*xiyf.^3; xi_yyyy=(12/a)*xiyf.^4;
  D1=(xi_y*ons).*D1xi; D2=((xi_y.^2)*ons).*D2xi + (xi_yy*ons).*D1xi;
  D3=((xi_y.^3)*ons).*D3xi+3*((xi_y.*xi_yy)*ons).*D2xi +(xi_yyy*ons).*D1xi;
  D4=((xi_y.^4)*ons).*D4xi+6*((xi_y.^2.*xi_yy)*ons).*D3xi + ...
      ((3*xi_yy.^2 + 4*xi_y.*xi_yyy)*ons).*D2xi + (xi_yyyy*ons).*D1xi;
else
  y=(xi+1)/2; D1 = 2*D1xi; D2 = 4*D2xi; D3 = 8*D3xi; D4 = 16*D4xi;
end
U = 1-y; Uy = -ones(N1,1); Uyy = zeros(N1,1);
ar2 =al^2*y; d =ar2+n^2; d2 =(l+1)*d-2*ar2; B21 =zeros(N1);
if (n==0),  B12 = B21;  A12 = B12;  A21 = B21;  B22 = D0; end;
if ((n==0) && (al~=0)), B11 = 8*D1 + 4*(y*ons).*D2-al*al*D0;
    A22 = Re*al*(U*ons).*D0 +1i*B11; A11 = Re*al*(U*ons).*B11 -Re*4*al ...
        *((y.*Uyy)*ons).*D0 +1i*( ((4*y).^2*ons).*D4 + 96*(y*ons).*D3 ...
        + ((96 - 8*ar2)*ons).*D2 - 16*al^2*D1 + al^4*D0 );
elseif ((n==0) && (al==0)), B11 =D0; A11 = 1i*4*(D1 + (y*ons).*D2);
    A22 = 1i*4*(2*D1 + (y*ons).*D2);
end
if (n~=0), g1 =(n^2+2*(l+1))*n^4 + ar2.*((3*n^2+4*(l+1))*n^2 ...
        +ar2.*(3*n^2+2*(l+1)+ ar2));%(g1-g15 are as defined in above paper)
    g2 =4*( (l+2)*n^6 + ar2.*((3*l+5)*n^4+ar2.*((3*l+4)*n^2+ar2*(l+1))));
    g3 =(l^4+8*l^3+26*l^2+32*l+13)*n^2+ar2.*(3*l^4+20*l^3+50*l^2+36*l+3 ...
        + ar2.*(3*l^2+10*l+7+ ar2)); g4=(l^3+6*l^2+11*l+6)*n^4+ar2.*( ...
        (3*l^3+15*l^2+21*l+5)*n^2+ar2.*(3*l^3+12*l^2+13*l+4 + ar2*(l+1)));
    g5 =8*((2*l^2+10*l+12)*n^6+ar2.*((5*l^2+22*l+21)*n^4+ar2.*( ...
        (3*l^2+12*l+9)*n^2+ar2.*(-n^2-ar2)))); g6 =d + 2*(l+1);
    g7=n^2*(l^2+4*l+7)+ar2.*(2*(l^2+3*l+2)+ar2); g8=n^2*(3*n^2+4*(l+1));
    g9=n^4*(2*(5*l+9)-n^2); g10=4*n^6*(l+3); g11=n^4*(n^2+2*l+18);
    g12=3*l^4+20*l^3+50*l^2+36*l+3; g13=n^2*(l^4-16*l^3-94*l^2-136*l-27);
    g14=8*n^4*(l^3+l^2-11*l-15); g15 = 16*n^6*(l^2+7*l+12);


    B11 = -al^2*(g1*ons).*D0 +(g2*ons).*D1 + 4*((y.*d.^3)*ons).*D2;
    B12 = -2*al*n*((d.^2)*ons).*D0; B22 = ((d.^2)*ons).*D0;
    A11 = Re*al*(U*ons).*B11 - Re*4*al*((y.*d.^3.*Uyy)*ons).*D0 ...
        - Re*8*al*n^2*((d.^2.*Uy)*ons).*D0 +1i*( ((al^4*g3)*ons).*D0 ...
        - 8*((al^2*g4)*ons).*D1 + (g5*ons).*D2 ...
        + 8*((y.*(g2+4*d.^3))*ons).*D3 + 16*(((y.*d).^2.*d)*ons).*D4 );
    A12 = Re*al*(U*ons).*B12 +1i*( -8*al^3*n*(d2*ons).*D0 ...
        -16*al*n*((d.*ar2)*ons).*D1 );
    A21 = Re*2*n*((Uy.*d.^2)*ons).*D0 +1i*( 2*al^3*n*(g6*ons).*D0 ...
        -4*al*n*((d2+(l+3)*d)*ons).*D1 -8*al*n*((y.*d)*ons).*D2 );
    A22 = Re*al*((U.*d.^2)*ons).*D0 +1i*(-al^2*(g7*ons).*D0 ...
        +4*((d.*((l+1)*d+n^2))*ons).*D1 + 4*((y.*d.^2)*ons).*D2);
end

% For n~=0, the system of phi and Omega are coupled.
if n~=0, % Enforcing no-slip and regularity Conditions:
    A = [[D0(1,:);D1(1,:);A11(3:N1,:)] [zeros(2,N1); A12(3:N1,:)]; ...
        [zeros(1,N1); A21(2:N1,:)]  [D0(1,:); A22(2:N1,:)]];
    B = [[zeros(2,N1); B11(3:N1,:)]  [zeros(2,N1); B12(3:N1,:)]; ...
        [zeros(1,N1); B21(2:N1,:)]  [zeros(1,N1); B22(2:N1,:)]];
    Breg_phi=-al^4*g8*D0(N1,:) +al^2*g9*D1(N1,:)+g10*D2(N1,:);
    Breg_omg=-4*al^3*n^3*D0(N1,:)-2*al*n^5*D1(N1,:); Areg_phi ...
        =Re*al*Breg_phi+Re*al^3*g11*D0(N1,:)-Re*4*al*l*n^6*D1(N1,:) ...
        +1i*( al^6*g12*D0(N1,:)+al^4*g13*D1(N1,:)-al^2*g14*D2(N1,:) ...
        +g15*D3(N1,:)); Areg_omg=Re*al*Breg_omg+Re*2*al^2*n^5*D0(N1,:) ...
        +1i*(-8*al^5*n*(l-1)*D0(N1,:) - 8*al^3*n^3*(l+3)*D1(N1,:));
    A(N,:) = [Areg_phi Areg_omg]; B(N,:) = [Breg_phi Breg_omg];
    %Gather no-slip conditions at top rows of the system [Aq =e Re Bq]:
    tmp1=A(3:N1,:); A(3,:) =A(N1+1,:); A=[A(1:3,:); tmp1; A(N1+2:end,:)];
    tmp1=B(3:N1,:); B(3,:) =B(N1+1,:); B=[B(1:3,:); tmp1; B(N1+2:end,:)];
    %Gather chosen 3 unknowns of q to be eliminated at first 3 locs. of q:
    tmp1 = A(:,3:N1); A(:,3) =A(:,N1+1); A = [A(:,1:3) tmp1 A(:,N1+2:end)];
    tmp1 = B(:,3:N1); B(:,3) =B(:,N1+1); B = [B(:,1:3) tmp1 B(:,N1+2:end)];
    %Remove singularities due to no-slip in [Aq = e Re Bq] by elimination:
    for kk = 1:3, vec = (A(kk,:)/A(kk,kk)); jj = kk+1:2*N1;
        A(jj,:) = A(jj,:)-A(jj,kk)*vec; B(jj,:) = B(jj,:)-B(jj,kk)*vec;
    end % Gauss elimination ended; Now solve the reduced system:
    L = (B(4:end,4:end)\A(4:end,4:end))/Re; [q,e] = eig(L); e =diag(e);
    ok=((imag(e)>-1*al)&(al>=1))|((imag(e)>-1)&(al<1));e=e(ok); q =q(:,ok);
    opts.tol=1e-21;opts.maxit=1e4;
    for ie=1:length(e), opts.v0=q(:,ie); e(ie)=eigs(L,1,e(ie),opts);  end
    [eimag,I]=sort(-imag(e));e=e(I);q=q(:,I);tmp1=q(1:N-1,:);
    tmp2=q(N:end,:);q_temp=-A(1:3,1:3)\ (A(1:3,4:end)*q);%Elim'ted coefs.
    q(1:N1,:) = [q_temp(1:2,:); tmp1];q(1+N1:2*N1,:) = [q_temp(3,:); tmp2];
elseif ((n==0) && (al~=0)) %Enforce Conditions for this and the other case:
    A11=[D0(1,:);D1(1,:);A11(3:N1,:)]; B11 = [zeros(2,N1); B11(3:N1,:)];
    A22 = [D0(1,:); A22(2:N1,:)]; B22 = [zeros(1,N1); B22(2:N1,:)];
    A11(N,:)=al*Re*(12*D2(N1,:)-(al^2+8)*D1(N1,:)+al^2*D0(N1,:))...
        +1i*(192*D3(N1,:)-24*al^2*D2(N1,:)+al^4*D1(N1,:));
    B11(N,:)=12*D2(N1,:)-al^2*D1(N1,:);
else A11=[D0(1,:);A11(2:N1,:)]; B11 = [zeros(1,N1); B11(2:N1,:)];
    A22 = [D0(1,:); A22(2:N1,:)]; B22 = [zeros(1,N1); B22(2:N1,:)];
end

% For n==0, the systems are decoupled. So solve for q1 = [phi or psi1] and
% q2 = [Omega or psi2] individually for accuracy and speed.
if n==0, M =2*(al~=0)+1*(al==0);%M= 2 or 1; M is No. of zero-rows in B_11
    for kk = 1:M, vec = (A11(kk,:)/A11(kk,kk));%Gauss elimination
        jj = kk+1:N1; A11(jj,:) = A11(jj,:)-A11(jj,kk)*vec;
        B11(jj,:) = B11(jj,:)-B11(jj,kk)*vec;
    end
    L1 = (B11(M+1:end,M+1:end)\A11(M+1:end,M+1:end))/Re; [q1,e1] = eig(L1);
    e1 =diag(e1);ok =((imag(e1)>-1*al)&(al>=1))|((imag(e1)>-1)&(al<1));
    e1 = e1(ok); q1 =q1(:,ok); opts.tol=1e-21;opts.maxit=1e4;
    for ie=1:length(e1),opts.v0=q1(:,ie);e1(ie)=eigs(L1,1,e1(ie),opts); end
    [e1imag,I]=sort(-imag(e1));e1=e1(I);q1=q1(:,I);
    q1=[-A11(1:M,1:M)\ (A11(1:M,M+1:end)*q1); q1];
    vec = (A22(1,:)/A22(1,1)); A22(2:N1,:) = A22(2:N1,:)-A22(2:N1,1)*vec;
    B22(2:N1,:) = B22(2:N1,:)-B22(2:N1,1)*vec;
    L2 = (B22(2:end,2:end)\A22(2:end,2:end))/Re; [q2,e2] = eig(L2);
    e2 =diag(e2); ok =((imag(e2)>-1*al)&(al>=1))|((imag(e2)>-1)&(al<1));
    e2 = e2(ok); q2 =q2(:,ok);
    for ie=1:length(e2),opts.v0=q2(:,ie);e2(ie)=eigs(L2,1,e2(ie),opts); end
    [e2imag,I]=sort(-imag(e2));e2=e2(I);q2=q2(:,I);
    q2=[-A22(1,1)\ (A22(1,2:end)*q2); q2];
    e = [e1; e2]; [eimag,I]=sort(-imag(e));e=e(I);
    q = [[q1; zeros(N1,size(q1,2))] [zeros(N1,size(q2,2)); q2]]; q=q(:,I);
end
\end{lstlisting}}
\bibliographystyle{elsarticle-num}
\bibliography{malik_skote}
\end{document}